\documentclass[preprint]{imsart}

\usepackage{amsthm,amsmath,amssymb,graphicx,enumitem,afterpage,bbm,comment}
\RequirePackage[numbers]{natbib}
\RequirePackage[colorlinks,
             linkcolor=ultramarine,
             linktoc=ultramarine,
             citecolor=ultramarine,
             urlcolor=ultramarine,
             filecolor=black]{hyperref}

\arxiv{1712.03522}

\startlocaldefs

\newtheorem{theorem}{Theorem}[section]

\newtheorem{lemma}[theorem]{Lemma}
\newtheorem{prop}[theorem]{Proposition}
\newtheorem{remark}[]{Remark}
\numberwithin{equation}{section}      
\numberwithin{remark}{section}      
\usepackage[math]{easyeqn}

\renewcommand{\(}{$\,}
\renewcommand{\)}{\,$}

\def\eqdef{\stackrel{\operatorname{def}}{=}}

\renewcommand{\tilde}[1]{\widetilde{#1}}

\renewcommand{\Gamma}{\varGamma}
\renewcommand{\Pi}{\varPi}
\renewcommand{\Sigma}{\varSigma}
\renewcommand{\Delta}{\varDelta}
\renewcommand{\Lambda}{\varLambda}
\renewcommand{\Psi}{\varPsi}
\renewcommand{\Phi}{\varPhi}
\renewcommand{\Theta}{\varTheta}
\renewcommand{\Omega}{\varOmega}
\renewcommand{\Xi}{\varXi}
\renewcommand{\Upsilon}{\varUpsilon}

\usepackage{color}

\definecolor{blue(pigment)}{rgb}{0.2, 0.2, 0.6}
\definecolor{ultramarine}{rgb}{0.07, 0.04, 0.56}
\definecolor{darkspringgreen}{rgb}{0.09, 0.45, 0.27}
\definecolor{hookersgreen}{rgb}{0.0, 0.44, 0.0}
\definecolor{plum(traditional)}{rgb}{0.56, 0.27, 0.52}
\definecolor{purple(html/css)}{rgb}{0.5, 0.0, 0.5}
\definecolor{magenta(dye)}{rgb}{0.79, 0.08, 0.48}

\newcommand{\data}{\boldsymbol{X}^n}
\newcommand{\subspaceset}{\mathcal{J}}
\newcommand{\Su}{\boldsymbol{\Sigma}}
\newcommand{\Se}{\widehat{\Su}}
\newcommand{\St}{{\Su^{*}}}
\newcommand{\Sp}{\widetilde{\Su}}

\newcommand{\Ou}{\boldsymbol{\Omega}}

\newcommand{\Oe}{\widehat{\Ou}}
\newcommand{\Ot}{\Ou^{*}}
\newcommand{\Op}{\widetilde{\Ou}}

\newcommand{\G}{\boldsymbol{G}}
\newcommand{\Gg}{\boldsymbol{\widetilde{G}}}
\newcommand{\Q}{\boldsymbol{Q}}
\newcommand{\W}{\boldsymbol{W}}
\newcommand{\Hh}{\boldsymbol{H}}
\newcommand{\V}{\boldsymbol{V}}
\newcommand{\U}{\boldsymbol{U}}
\newcommand{\Aa}{\boldsymbol{A}}
\newcommand{\Bb}{\boldsymbol{B}}
\newcommand{\D}{\boldsymbol{D}}
\newcommand{\Ee}{\boldsymbol{E}}

\newcommand{\fl}{\rho}
\newcommand{\A}{\mathcal{A}}
\newcommand{\vic}{\mathcal{B}}
\newcommand{\res}{\varepsilon}

\def\Proj{\boldsymbol{P}}
\def\Projs{\Proj^{*}}
\def\Projh{\widehat{\Proj}}

\newcommand{\Pt}{\Projs_{\subspaceset}}
\newcommand{\Ptr}{\Projs_{s}}
\newcommand{\Pe}{\Projh_{\subspaceset}}

\newcommand{\Per}{\Projh_{s}}
\newcommand{\Pu}{\Proj_{\subspaceset}}
\newcommand{\Pur}{\Proj_{s}}

\def\Gammas{\Gamma^*}

\def\IndS{\mathcal{I}_{\subspaceset}}

\newcommand{\I}{\boldsymbol{I}}
\def\cond{\, \big| \,}
\def\Cond{\, \bigg| \,}

\def\prior{\Pi}
\newcommand{\Ppost}[1]{\prior\left( #1 \cond \data \right)}
\newcommand{\Ppostbig}[1]{\prior\left( #1 \Cond \data \right)}

\newcommand{\PpostW}[1]{\prior^{\Wsht}\left( #1 \cond \data \right)}
\newcommand{\PpostbigW}[1]{\prior^{\Wsht}\left( #1 \Cond \data \right)}
\newcommand{\Pro}{\mathbb{P}}
\newcommand{\R}{\mathbbm{R}}
\newcommand{\SPD}{\mathbbm{S}^{p}_{+}}
\newcommand{\E}{\mathbbm{E}}
\newcommand{\Var}{\operatorname{\mathbbm{V}ar}}
\newcommand{\Phib}{\boldsymbol{\Phi}}
\newcommand{\Psib}{\boldsymbol{\Psi}}
\newcommand{\Xiu}{\boldsymbol{\Xi}}
\newcommand{\Xit}{\boldsymbol{\Xi}^*}
\newcommand{\Xie}{\widehat{\boldsymbol{\Xi}}}

\newcommand{\sumn}{\sum\limits_{\textit{j}=1}^\textit{n}}

\newcommand{\ND}{\mathcal{N}}
\newcommand{\rank}{\operatorname{r}}
\newcommand{\reff}{\tilde{\operatorname{r}}}
\newcommand{\diag}{\operatorname{diag}}

\def\eqdist{\stackrel{d}{=}}
\def\Wsht{\mathcal{W}}
\def\IWsht{\mathcal{W}^{-1}}

\def\hdelta{\widehat{\delta}}
\def\wdelta{\widetilde{\delta}}

\def\Fr{2}
\def\T{\top}
\def\Tr{\operatorname{Tr}}

\def\CONST{\mathtt{C} \hspace{0.1em}}

\def\Deltas{\Delta^*}

\def\ms{m^{*}}
\def\gs{g^{*}}

\def\us{u^{*}}
\def\uh{\widehat{u}}

\def\mus{\mu^{*}}
\def\muh{\widehat{\mu}}
\def\sigmas{\sigma^{*}}
\def\sigmah{\widehat{\sigma}}

\endlocaldefs

\begin{document}

\begin{frontmatter}

\title{Finite sample Bernstein-von Mises theorems for functionals and spectral projectors of the covariance matrix}

\runtitle{BvM for functionals and projectors of covariance matrix}

\begin{aug}

\author{\fnms{Igor} \snm{Silin}\thanksref{a,b}\corref{}\ead[label=e1]{isilin@princeton.edu}}

\address[a]{Department of Operations Research and Financial Engineering, Princeton University, Princeton, NJ 08540, USA. \printead{e1}}
\address[b]{Institute for Information Transmission Problems RAS, Moscow, Russia}

\runauthor{I. Silin}

\affiliation{Princeton University}

\end{aug}

\begin{abstract}
    We demonstrate that a prior influence on the posterior distribution of the covariance matrix vanishes as sample size grows. 
The assumptions on priors are explicit and mild.
	The results are valid for finite sample and admit the dimension $p$ growing with sample size  $n$. 
	We exploit the described fact to derive the finite sample Bernstein-von Mises theorem for functionals of covariance matrix (e.g. eigenvalues) and to find the posterior distribution of the difference between spectral projector and empirical spectral projector. 
	This can be useful for constructing sharp confidence sets for the true value of the functional or for the true spectral projector.
\end{abstract}




\end{frontmatter}

    \section{Introduction} \label{Section: Introduction}
		The Bernstein-von Mises (BvM) phenomenon states some pivotal behaviour of the posterior distribution.
It specifies conditions on a prior, under which the influence of the prior vanishes as the number of observations grows and the posterior is asymptotically Gaussian.
The main application of BvM is the usage of Bayesian credible sets as frequentist confidence sets.
It helps in situations when the frequentist uncertainty quantification does not allow to build confidence sets directly due to unknown parameters of the asymptotic distribution.

The classical books \cite{LeCam,Van der Vaart} contain chapters on BvM for standard parametric setup. 
More general semiparametric models were studied in \cite{Bickel_BvM}.
In modern statistics main focus is on the growing parameter dimension, so the classical results should be reconsidered; 
see, e.g. \cite{Ghosal,Johnstone_BvM,Spokoiny_BvM} for some examples in high dimensions.
Moreover, modern statisticians are focused on models with samples of limited size, however, only a few finite sample BvM results are available, e.g. \cite{Panov}.
We also mention the BvM for linear functionals of the density derived in \cite{Rivoirard}, general theory for smooth functionals of the target parameter presented in \cite{Castillo_BvMFSFISM}, as well as \cite{Castillo_NBvMTIGWN,Castillo_BvMPNBP} and some recent works on BvM for statistical inverse problems \cite{Nickl_1,Nickl_2,Giordano}, among other important results in the field. The frequentist coverage properties of Bayesian credible sets are studied in \cite{Szabo2,Szabo}.

This paper aims at deriving similar results for the following specific model. 
Let the data $ \data = (X_1, \ldots, X_{n}) $ be independent identically distributed 
zero-mean random vectors in $ \R^p $.
Its covariance matrix is given by
\begin{equation}
	\begin{aligned}
		\St \eqdef \E \left(X_j X_j^{\T}\right).
	\nonumber
	\end{aligned}
\end{equation}
The natural estimate of the unknown true covariance is the sample covariance matrix, defined as
\begin{equation}
    \begin{aligned}
    	\Se \eqdef \frac{1}{n} \sumn X_j X_j^{\T}.
    \nonumber
    \end{aligned}
\end{equation}
The quality of estimation in spectral norm $ \| \Se - \St \|_{\infty} $ arises in numerous problems and is well-examined; 
see, for instance, \cite{Koltchinskii_CIAMBFSCO, Rigollet, Tropp, Vershynin_ITTNAAORM, Adamczak}. 
Functionals and spectral projectors of covariance matrix also appear in applications frequently, so this model is of special interest.

In this work we show that the posterior distribution of the covariance matrix $ \Su $ stays approximately the same for different choices of prior distribution as soon as we have enough observations. 
In particular, we demonstrate that the posterior computed from an arbitrary prior deviates not significantly from the posterior that comes from one special class of priors -- the Inverse Wishart prior, which is the conjugate to the multivariate Gaussian distribution.
We do not impose any conditions on a prior; the error of approximation (which we measure by the total variation distance) of the corresponding posterior by the Inverse Wishart posterior is described in terms of two crucial concepts. 
One of these concepts is the \textit{posterior contraction}, and the other is the \textit{flatness} of the prior. 
So, our result makes sense if the prior is such that:
\begin{itemize}
	\item the posterior concentrates in a relatively small vicinity of the true covariance $ \St $;
	\item the prior is ``flat'' enough in this vicinity, that is, can be well-approximated by a constant.
\end{itemize}

The described ``posterior independence'' enables the following strategy that allows to reduce the complexity of a problem at hand drastically.
Instead of working with some complicated prior, one can consider the Inverse Wishart prior. 
Since this prior is the conjugate to the multivariate Gaussian distribution, the posterior is again the Inverse Wishart, so we can study it directly.
Moreover, in a wide range of situations nice properties of the Inverse Wishart distribution simplify the analysis significantly.
 
We apply the proposed strategy to the following objects that are of crucial importance in modern applications.
First, we derive the BvM theorem for approximately linear functionals of covariance matrix. 
The main focus here is on eigenvalues of the covariance matrix, which are extensively studied object, see \cite{MP,El Karoui,Johnstone_TW} and references therein. 
The asymptotic normality of the posterior measure of approximately linear functionals was already shown in \cite{Zhou_BvMTFFOCM}, which is based on the technique developed in \cite{Castillo_BvMFSFISM}.
However, not only does their result require sample size to grow to infinity, but also imposes some nontrivial condition on a prior instead of our simple ``flatness'' assumption. 

The second object under consideration is a spectral projector of the covariance matrix.
Let $ \Pt $ be the projector onto a subset $ \subspaceset $ of eigenspaces of $ \St $.
Its sample version is given by $ \Pe $ corresponding to the sample covariance $ \Se $. 
For some recent results on the distribution of $ \| \Pe - \Pt \|_{\Fr}^2 $ we refer to \cite{Koltchinskii_NAACOSPOSC,Lounici_new,Naumov_BCSFSPOSC,Silin}.
These objects are closely related to the Principal Component Analysis (PCA), probably the most famous dimension reduction method. PCA-based methods are widely used in finance \cite{Pasini}, as well as in other fields, see \cite{Fan_PCA} for an overview of statistical problems related to PCA and references therein.
Recent developments in theoretical guarantees for sparse PCA in high dimensions engender attention to such methods, see \cite{Johnstone,Birnbaum,Berthet_ODSPCHD,Cai_SPCAORAE,Gao_ROPCSPCA}. 
Our approach is applied to the squared Frobenius distance $ \| \Pu - \Pe \|_{\Fr}^2 $ between the projector of $ \Su $ and the projector of $ \Se  $.
It is interesting to mention that while the posterior distribution of a functional is approximated by Gaussian, which is usual for BvM,
for the case of spectral projectors the limiting distribution is the distribution of Gaussian quadratic form.

It is compelling that the presented technique does not rely on Gaussianity of the data.
Even though we work with the Gaussian likelihood, we allow \textit{model misspecification} and formulate our results for \textit{pseudo-posterior}.
As to the distribution of the data, we require only one property: the concentration of the sample covariance $ \Se $ around the true covariance $ \St $. 

Also, even though the assumptions that the eigenvalues of $ \St $ are bounded from above and separated from zero, or that the spectral gaps (differences between consequent eigenvalues) are separated from zero are common in literature, we avoid them.
Our results admit growing spectral norm $ \| \St \|_{\infty} $ and vanishing smallest eigenvalue and spectral gaps.
The provided error bounds are explicit and allow to track what regimes of $ \St $ still ensure convergence to the limiting distribution.

The main contributions of this paper are as follows.

\begin{itemize}
\item
We establish a finite sample result stating that the prior influence on the posterior distribution of covariance matrix disappears as the sample size grows.
The assumptions on a prior are mild and easy to verify. 
The distribution of the data can also be quite general. 
\item 
We propose a novel strategy for analyzing the posterior distribution for an arbitrary prior. 
The strategy includes: first, approximation of our posterior at hand by the posterior based on the conjugate prior, and second, study of the latter posterior which has nice properties.
\item
The described strategy is applied to derive finite sample BvM theorems for functionals and spectral projectors of the covariance matrix. In case of linear functionals, the obtained error rates are dimension-free for the conjugate prior and for the special class of ``rank-adjusted'' non-conjugate priors that we present.
\end{itemize}

The rest of the paper is structured as follows. 
The rest of this section, that is, Subsection~\ref{Subsection: Notations}, contains some notations.
Section~\ref{Subsection: Setup and problem} explains our setup.
Bayesian framework is described in Section~\ref{Subsection: Bayesian framework and main result}. 
Section~\ref{Subsection: Application} is dedicated to the applications of the proposed strategy: functionals of covariance matrix are studied in Subsection~\ref{Subsection: functionals}, and spectral projectors are considered in Subsection~\ref{Subsection: projectors}. We devote Section~\ref{Section: Discussion} to the discussion and comparison of the obtained rates.
The proofs of the main theorems are collected in Section~\ref{Section: Proofs}.
Appendix~\ref{Appendix: A} and Appendix~\ref{Appendix: B} gather some auxiliary results from the literature and the rest of the proofs, respectively.
		\subsection{Notations} \label{Subsection: Notations}
			We will use the following notations throughout the paper.
The space of real-valued $d_1 \times d_2$ matrices is denoted by $\R^{d_1 \times d_2}$, 
while $\mathbbm{S}^{d}_{+}$ means the set of positive-semidefinite matrices of size $d \times d$. 
We write $ \I_d $ for the identity matrix of size $ d \times d $,
$ \rank(A) $ and $ \Tr(B) $ stand for the \emph{rank} of a matrix $ A $ and the \emph{trace} of a square matrix $ B $.
Further, we stick to the notations of the Schatten norm, so that $ \| A \|_{\infty} $ stands for the \emph{spectral norm} of a matrix $A$,
while $\| A \|_1$ means the \emph{nuclear norm}.
The \emph{Frobenius scalar product} of two matrices $ A $ and $ B $ of the same size is $ \langle A, B \rangle_2 \eqdef \Tr(A^{\T}B) $,  
while the \emph{Frobenius norm} is denoted by $ \| A \|_{2} $. 
When applied to a vector, $ \| \cdot \| $ means just its \emph{Euclidean norm}.
The \emph{effective rank} of a square matrix $ B $ is defined by $ \reff(B) \eqdef \frac{\Tr(B)}{\| B\|_{\infty}} $. The \emph{condition number} of a square invertible matrix $B$ is $\kappa(B) \eqdef \| B \|_{\infty} \| B^{-1} \|_{\infty}$.

The relation $ a \lesssim b $ means that there exists an absolute constant $ C $, different from place to place, such that $ a \leq Cb $, while
$ a \asymp b $ means that $ a \lesssim b $ and $ b \lesssim a $.
By $ a\lor b $ and $ a\land b $ we mean maximum and minimum of $ a $ and $ b $, respectively.
In the sequel we will often be considering intersections of events of probability greater than $ 1-1/n $. 
Without loss of generality, we will write that probability measure of such an intersection is $ 1 - 1/n $, since it can be easily achieved by adjusting constants.
Finally, we write $\eta \eqdist \xi$ when the distributions of the random variables $\eta$ and $\xi$ coincide, and $ \eta_{n} = o_{\Pro}(1) $ when the sequence $\{\eta_{n}\}_{n=1}^{\infty}$ of random variables converges to $0$ in probability. 
		
	\section{Setup and main result} \label{Section: Procedure and main results}
		This section explains our setup and states the main result.
		\subsection{Setup} \label{Subsection: Setup and problem}
			Recall that we consider the model where the observations $ \data = (X_1, \ldots, X_{n}) $ are independent identically distributed 
zero-mean random vectors in $ \R^p $ with the true covariance matrix
\begin{equation}
	\begin{aligned}
		\St \eqdef \E \left(X_j X_j^{\T}\right)
	\nonumber
	\end{aligned}
\end{equation}
and the sample covariance matrix
\begin{equation}
    \begin{aligned}
    	\Se \eqdef \frac{1}{n} \sumn X_j X_j^{\T}.
    \nonumber
    \end{aligned}
\end{equation}
Assume without loss of generality that $ \St \in \SPD $ is full rank and invertible (otherwise, using $\Se$, one can easily pass to a subspace of $\R^p$ where the true covariance matrix for the transformed data will be of full rank). 

We do not need the assumption on Gaussianity of the data.
The only condition that our main result require from the underlying distribution of the data $ \data = (X_1, \ldots, X_{n}) $ is the concentration of the sample covariance matrix $ \Se $ around the true covariance $ \St $:
\begin{equation}
    \begin{aligned}
        \| \Se - \St\|_{\infty} \leq \hdelta_{n,p} \|\St\|_{\infty}
    \label{Eq: hat_delta_n}
    \end{aligned}
\end{equation}
with probability $ 1 - {1}/{n} $. 
The bound $ \hdelta_{n,p} $ from the condition depends on the sample size $n$, the dimension $p$ and potentially $\St$, but we omit this dependence to keep the notation light. Clearly, $ \hdelta_{n,p} $  varies for different distributions of the data, but it allows to work with much wider classes of probability measures than just Gaussian or sub-Gaussian. While for the Gaussian case one may take
\begin{equation}
    \begin{aligned}
        \hdelta_{n,p} \;\asymp \;  \sqrt{\frac{\reff(\St) + \log(n)}{n}} ,
        \nonumber
    \end{aligned}
\end{equation}
several more examples of possible distributions and the corresponding $ \hdelta_{n,p} $ for them are provided in Theorem \ref{Th: Covariance concentration} from Appendix \ref{Appendix: A}. 
Throughout the rest of the paper we assume that the data satisfy condition (\ref{Eq: hat_delta_n}).
		\subsection{Bayesian framework and main result} \label{Subsection: Bayesian framework and main result}
			In Bayesian framework one imposes a prior distribution $ \Pi $ on the covariance matrix $ \Su $. 
Even though our data are not Gaussian, we can consider the Gaussian log-likelihood:
\begin{equation}
    \begin{aligned}
    	l_{n}(\Su) = -\frac{n}{2}\log{\det(\Su)} - \frac{n}{2} \Tr(\Su^{-1}\Se) - \frac{np}{2} \log{(2\pi)}.
    \nonumber
    \end{aligned}
\end{equation}
The posterior measure of a set $ \A \subset \SPD $ can be expressed as
\begin{equation}
    \begin{aligned}
    	\Ppost{\A} = \frac{\int_{\A} {\exp{\left(l_{n}(\Su) \right)} \;d\Pi(\Su)}}{\int_{\SPD}{\exp{\left(l_{n}(\Su)\right)} \;d\Pi(\Su)}}.
    \nonumber
    \end{aligned}
\end{equation}
As the Gaussian log-likelihood $ l_{n}(\Su) $ does not necessarily correspond to the true distribution of our data, 
we call the random measure $ \Ppost{\;\cdot\;} $ a \textit{pseudo-posterior}.
Once the prior is fixed, we can easily sample matrices $ \Su $ from this pseudo-posterior distribution. In what follows, we drop the prefix pseudo- and call $ \Ppost{\;\cdot\;} $ just posterior.

Our technique relies on two crucial concepts. 
The first one is \textit{posterior contraction}. Define $ \delta-$vicinity of $ \St $ as:
\begin{equation}
    \begin{aligned}
    \vic(\delta) 
     \eqdef 
    \left\{ \Su \in \SPD: \; \| \Su - \St \|_{\infty} \leq \delta \, \| \St\|_{\infty} \right\}.
    \label{Vicinity}
    \end{aligned}
\end{equation}
Then we can find a radius $ \delta_{n,p} $ such that the following posterior contraction condition is fulfilled:
\begin{equation}
    \begin{aligned}
    \Ppost{\vic(\delta_{n,p})} 
    \geq 
    1 - \frac{1}{n}
    \label{Def: Contraction}
    \end{aligned}
\end{equation}
with probability $ 1 - {1}/{n}$. Proving the posterior contraction and finding the radius $\delta_{n,p}$ for an arbitrary prior is a separate nontrivial problem of independent interest, which lies beyond the scope of this work.

The second concept that we introduce is \textit{``flatness''} of the prior, defined as
\begin{equation}
    \begin{aligned}
    \fl(\delta) 
    \eqdef 
    \sup\limits_{\Su \in \vic(\delta)} \left| \frac{\Pi(\Su)}{\Pi(\St)} - 1\right|.
    \label{Def: flatness}
    \end{aligned}
\end{equation}
It is typically easy to bound when we have the closed-form expression for the density of the prior.

We will be actively using the conjugate prior $\Pi^{\Wsht}$ to the multivariate Gaussian distribution, that is, the Inverse Wishart distribution $ \IWsht_p(\G,\, p+b-1) $ with parameters $ \G\in\SPD,\, \G \succ 0 $ and $b\in\R, \,b > 0$. 
Its density is given by
    \begin{equation}
        \begin{aligned}
            \frac{d\Pi^{\Wsht}(\Su)}{d\Su} \propto \exp{ \left( -\frac{2p+b}{2}\log \det(\Su) - \frac{1}{2} \Tr(\G\Su^{-1})\right)}.
        \nonumber
        \end{aligned}
    \end{equation}
\begin{remark} 
Since $b$ is under our control, it makes sense to assume $b \asymp 1$ only to avoid annoying terms in the error bounds; similarly, $\G$ is chosen by us, so $\| \G \|_{\infty}$ and $\| \G \|_1$ can be made arbitrarily small. Hence, the terms including $b$ or $\G$ in the sequel may always be considered negligible.
\end{remark}
Useful basic properties of the Inverse Wishart distribution are stated in Lemma~\ref{Lemma: Properties} relegated to Appendix~\ref{Appendix: A}.
The following two propositions state the posterior contraction and a bound for the flatness of this special prior, respectively.
\begin{prop} \label{Contraction}
    Let $\delta_{n,p}^{\Wsht}$ be the posterior contraction radius (\ref{Def: Contraction}) of the Inverse Wishart prior $\Pi^{\Wsht}$. Assume $p/n$ is smaller than some implicit constant. Then the following bound holds:
    \begin{equation}
        \begin{aligned}
        \delta_{n,p}^{\Wsht} 
         \lesssim 
        \sqrt{\frac{\log(n)+p}{n}} + \hdelta_{n,p} + \frac{\|\G\|_{\infty}}{n\,\|\St\|_{\infty}} \,. 
        \nonumber
        \end{aligned}
    \end{equation} 
\end{prop}
\begin{prop} \label{Flatness}
    Let $ \fl^{\Wsht}(\delta) $ be the flatness (\ref{Def: flatness}) of the Inverse Wishart prior $\Pi^{\Wsht}$.
    For $\delta\leq1/(2\,p\,\kappa(\St) )$ the following bound holds:
    \begin{EQA}
        \fl^{\Wsht}(\delta) 
        & \lesssim &
        \delta \,\kappa(\St) \,\left( p^2 + \| \G \|_1 \| \St^{-1} \|_{\infty} \right).
    \end{EQA}
\end{prop}
\noindent The proofs are postponed to Appendix~\ref{Appendix: B}.

The first main result of this paper shows that the posterior distribution is approximately the same for sufficiently flat priors satisfying the posterior contraction condition. In particular, we bound the total variation distance between the posterior distributions computed from an arbitrary non-conjugate prior and the conjugate prior.
\begin{theorem} \label{Theorem: posterior independence}
    Assume the distribution of the data $ \data = (X_1, \ldots, X_{n}) $ fulfills the sample covariance concentration property (\ref{Eq: hat_delta_n}).
    Consider an arbitrary prior $\Pi$ and the prior $ \Pi^{\Wsht} $ given by the Inverse Wishart distribution $ \IWsht_p(\G, \,p+b-1) $.
    Define 
    \begin{EQA}
        \overline{\delta}_{n,p} 
        & \eqdef &
        \delta_{n,p} \, \lor \, \delta_{n,p}^{\Wsht}
    \end{EQA}
    with $ \delta_{n,p} $ and $ \delta_{n,p}^{\Wsht} $ from (\ref{Def: Contraction}) for $\Pi$ and $ \Pi^{\Wsht} $, respectively.
    Then the following holds with probability $ 1 - 1/n $:
    \begin{EQA} 
            \sup\limits_{\A \subset \SPD} \left| \Ppost{\A} - \PpostW{\A} \right| 
            & \lesssim &
            \Diamond^{\ast}_{\Pi} \, ,
    \end{EQA}
    where
    \begin{equation}
        \begin{aligned}
        \Diamond^{\ast}_{\Pi} 
        \eqdef 
        \fl( \overline{\delta}_{n,p} ) + \fl^{\Wsht}( \overline{\delta}_{n,p} ) + \frac{1}{n} \, .
        \label{Def: Diamond_ast} 
        \end{aligned}
    \end{equation}
\end{theorem}
\begin{remark} \label{Remark BvM}
     Proposition~\ref{Contraction} tells us that the usual rate for posterior contraction radius is $\sqrt{(\log(n)+p)/n} + \hdelta_{n,p}$. Proposition~\ref{Flatness} helps to understand that we can expect linear behaviour of the flatness $\fl(\delta)$ on $ \delta$. Therefore, if the prior behaves properly (i.e. $\delta_{n,p}$ of the same order as $\delta_{n,p}^{\Wsht}$ and $\rho(\delta)$ can be bounded, say, by a linear function in $\delta$), we can expect 
    \begin{EQA}
        && \Diamond^{\ast}_{\Pi} \asymp \CONST(\St) \, \left( \sqrt{\frac{p^4\log(n)+p^5}{n}} + p^2\,\hdelta_{n,p} \right).
    \end{EQA}
    Our result makes sense if $p^5/n \ll 1$, which is, undoubtedly, very strict requirement; however, this is not due to inaccurate bounds in the proof. In general, the technique that we propose doesn't allow to state a better bound, probably because it is too general; more discussion on that is available in Section~\ref{Section: Discussion}. At the same time, in some simple cases we will be able to get a dimension-free bound for special class of non-conjugate priors constructed in Subsection~\ref{Rank-adjusted Prior}.
\end{remark}
	\section{Examples} \label{Subsection: Application}
		Before considering particular examples, let us introduce some additional notations concerning the true covariance $\St$ .

Let $\sigmas_1 \geq \ldots \geq \sigmas_p$ be the ordered eigenvalues of $\St$ .
Suppose among them there are $q$ distinct eigenvalues $ \mus_1 > \ldots > \mus_q $.
Introduce groups of indices $\Deltas_{s} = \{ j: \mus_{s} = \sigmas_j\}$ and denote by $\ms_{s}$ the multiplicity $|\Deltas_{s}|$  for all $s \in \{1,\ldots, q\}$ . 
The corresponding eigenvectors are denoted as $\us_1,\ldots, \us_p$ .
We will use the projector on the $s$-th eigenspace of dimension $\ms_{s}$ :
\begin{EQA}
    \Ptr 
    & = &
    \sum\limits_{j\in\Deltas_{s}} \us_j {\us_j}^{\T}
\end{EQA}
and the eigendecomposition
\begin{EQA}  
    \St 
    & = &
    \sum\limits_{j=1}^p \sigmas_j \us_j {\us_j}^{\T} = \sum\limits_{s=1}^q \mus_{s} \left( \sum\limits_{j\in\Deltas_{s}} \us_j {\us_j}^{\T} \right) =
        \sum\limits_{s=1}^q \mus_{s} \Ptr.
\end{EQA}
We also introduce the spectral gaps $\gs_s$ :
\begin{equation}
    \begin{aligned}
        \gs_s = \begin{cases}
         \mus_1 - \mus_{2}, & s = 1,\\
         (\mus_{s-1} - \mus_{s}) \; \land \; (\mus_s - \mus_{s+1}),\;\;& s  \in \{2,\ldots, q-1\},\\
         \mus_{q-1} - \mus_{q}, & s = q.
        \end{cases}     
    \nonumber
    \end{aligned}
\end{equation}  

Similarly, suppose that $\Se$  has $p$  (distinct with probability $1$) eigenvalues  $\sigmah_1 > \ldots > \sigmah_p$.
The corresponding eigenvectors are denoted as $\uh_1,\ldots, \uh_p$ .
Additionally, suppose that \\ $\| \Se - \St\|_{\infty} \leq \frac{1}{4} \min\limits_{s\in\overline{1,q}} \gs_s$.
Then, as shown in \cite{Koltchinskii_AACBFBFOSPOSC}, we can identify clusters of the eigenvalues of $\Se$ corresponding to each eigenvalue of $\St$ 
and therefore determine $\Deltas_s$ and $\ms_{s}$ for all $s \in \{1,\ldots, q\}$, so further we assume they are known.

		\subsection{Functionals of covariance matrix} \label{Subsection: functionals}
			\subsubsection{General case: conjugate and non-conjugate priors}
Represent a linear functional $ \phi(\Sp) $ as
\begin{equation}
    \begin{aligned}
    \phi(\Sp) - \phi(\St) 
    = 
    \Tr[\Phib (\Sp-\St)] + \res(\Sp,\St),
    \label{Linear expansion}
    \end{aligned}
\end{equation}
where we suppose there exists $ \Phib \in\SPD $ such that the residual $ \res(\Sp,\St) $ is bounded in the following way:
\begin{equation}
    \begin{aligned}
    | \res(\Sp,\St) |
     \leq  
    \CONST_{\phi}(\St) \, \| \Sp - \St\|_{\infty}^2,
    \label{Bound res}
    \end{aligned}
\end{equation}
with the constant $\CONST_{\phi}(\St)$ different for different functionals and depending only on $\St$.
So, in  some sense $ \phi(\cdot) $ is ``\textit{approximately linear}'' functional. Throughout the rest of the paper we denote
\begin{equation}
    \begin{aligned}
        r \eqdef \rank(\Phib)
    \nonumber
    \end{aligned}
\end{equation}
to be the rank of $\Phib$. It will play crucial role in our bounds.

Paper \cite{Zhou_BvMTFFOCM} provides several functionals that can be put in this context. The simplest of them, such as entries, quadratic forms, trace, are linear functionals with zero residual $\varepsilon(\,\cdot\,,\,\cdot\,)$, while functionals like log-determinant and eigenvalues are approximately linear. Let us briefly take a closer look on eigenvalues of the covariance matrix.
Define the functional as
\begin{EQA}
    \phi(\Sp) 
    & = &
    \frac{1}{\ms_s} \sum\limits_{j \in \Deltas_s} \widetilde{\sigma}_j.
\end{EQA}
Clearly, $ \phi(\St) = \mus_s $ and its natural estimate is 
\begin{EQA}
    \phi(\Se) = \muh_s 
    & \eqdef &
    \frac{1}{\ms_s} \sum\limits_{j \in \Deltas_s} \sigmah_j. 
\end{EQA}
The next proposition claims that the functional defined in such a way is approximately linear functional of covariance matrix.
\begin{prop} 
            Assume $\| \Sp - \St\| \leq \frac{\gs_s}{2e^2}$.
            Then the following bound  for first-order approximation takes place:
            \begin{equation}
                \begin{aligned}
                    \left| \widetilde{\mu}_s - \mus_s - \Tr \left[ \frac{\Ptr}{\ms_s}\;(\Sp - \St) \right] \right| \leq \frac{2e^2}{\gs_s} \, \| \Sp - \St \|^2,
                \nonumber
                \end{aligned}
            \end{equation}
            or, in other terms, introducing $\Phib = \frac{\Ptr}{\ms_s}$,
            \begin{equation}
                \begin{aligned}
                    | \res(\Sp, \St) | 
                    = 
                    \left| \phi(\Sp) - \phi(\St) - \Tr \left[ \Phib\,(\Sp - \St)  \right] \right| 
                    \leq 
                    \frac{2e^2}{\gs_s} \, \| \Sp - \St \|^2.
                \nonumber
                \end{aligned}
            \end{equation}
\end{prop}
So, for eigenvalues the assumption (\ref{Linear expansion}), (\ref{Bound res}) is fulfilled with $ \CONST_{\phi}(\St) = {2e^2}/{\gs_s} $.
We omit the proof of this statement and refer to \cite{Kato} for the details on perturbation theory for eigenvalues.
Let us continue with arbitrary functional $ \phi(\cdot) $ satisfying (\ref{Linear expansion}), (\ref{Bound res}), but keeping this eigenvalue example in mind.

The goal is to prove the finite sample BvM theorem  for a functional of covariance matrix by applying our general strategy described in the previous section.
We first need the result for the Inverse Wishart prior.

To be able to derive dimension-free result, we need to redefine the assumption of sample covariance concentration a little bit. Consider the matrix $\St^{1/2} \Phib \St^{1/2}$ of rank exactly $r$ and its spectral decomposition:
\begin{equation}
    \begin{aligned}
        \St^{1/2} \Phib \St^{1/2} = \V \D \V^{\top} \;\;\text{with}\;\; \V\in\R^{p\times r},\, \V^{\top} \V = \I_r \text{ and }\D\in\R^{r\times r} \text{ diagonal}.
        \nonumber
    \end{aligned}
\end{equation}
Define the quantity $\wdelta_{n,r}$ such that
\begin{equation}
    \begin{aligned}
        \| \V^{\top}\St^{-1/2} (\Se-\St) \St^{-1/2} \V \|_{\infty} \leq \wdelta_{n,r}
        \label{ConcentrationCondition2}
    \end{aligned}
\end{equation}
with probability $1-1/n$. Clearly, $\wdelta_{n,r}$ is the analogue of $\hdelta_{n,r}$ from (\ref{Eq: hat_delta_n}).
\begin{remark}
    We always have trivial inequality
    \begin{equation}
        \begin{aligned}
            \wdelta_{n,r} \leq \kappa(\St)\; \hdelta_{n,p}.
            \nonumber
        \end{aligned}
    \end{equation}
    However, the advantage of $\wdelta_{n,r}$ is that, unlike $\hdelta_{n,p}$, it doesn't depend on the dimension $p$ of the full space, since it describes the concentration of the sample covariance of $r$-dimensional random vectors $\V^{\top}\St^{-1/2}X_j,\;j=1,\ldots,n$. For instance, if our data are Gaussian, we always have
    \begin{equation}
        \begin{aligned}
            \wdelta_{n,r} \asymp \sqrt{\frac{r+\log(n)}{n}},
            \nonumber
        \end{aligned}
    \end{equation}
    while $\hdelta_{n,p}$ may be as large as
    \begin{equation}
        \begin{aligned}
            \hdelta_{n,p} \asymp \sqrt{\frac{p+\log(n)}{n}}.
            \nonumber
        \end{aligned}
    \end{equation}
\end{remark}
Now we are ready to state the functional Bernstein-von Mises theorem for conjugate prior. 
\begin{theorem} \label{Theorem: fBvM Wishart}
    Assume the distribution of the data $\data = (X_1, \ldots, X_{n})$ fulfills the sample covariance concentration properties (\ref{Eq: hat_delta_n}) and (\ref{ConcentrationCondition2}). 
    Suppose the functional $\phi$ satisfies (\ref{Linear expansion}), (\ref{Bound res}).
    Consider the prior $\Pi^{\Wsht}$ given by the Inverse Wishart distribution $\IWsht_p(\G, \; p+b-1)$.
    Let $\zeta \sim \ND(0,1)$.
    Then with probability $1 - 1/n$
    \begin{EQA}
        \sup\limits_{x\in\R} \left| \PpostbigW{\frac{\sqrt{n} \; \left( \phi(\Su) - \phi(\Se) \right)}{\sqrt{2} \| \St^{1/2} \Phib \St^{1/2} \|_2} \leq x} - 
        \Pro(\zeta \leq x)\right|
        & \lesssim &
        \Diamond_{\phi}, 
    \end{EQA}
    where
    \begin{equation}
        \begin{aligned}
        \Diamond_{\phi}
         \eqdef 
        \Diamond_{NL} + \Diamond_{GA}
        \label{Eq: Diamond_1}
        \end{aligned}
    \end{equation}
    and
    \begin{EQA}
        \Diamond_{NL} 
        & \eqdef &
        \sqrt{n} \left( \delta_{n,p}^{\Wsht} + \hdelta_{n,p} \right)^2 \,\frac{\CONST_{\phi}(\St) \, \| \St \|_{\infty}^2 }{\| {\St}^{1/2} \Phib {\St}^{1/2} \|_{\Fr}}\,,\\
        \Diamond_{GA}
        & \eqdef &
        \sqrt{r^2 + r\log(n)} \;\wdelta_{n,r} + \sqrt{\frac{r^3 + r^2 \log(n)}{n}} + \sqrt{\frac{r}{n}}\,\| \St^{-1/2}\G\St^{-1/2} \|_{\infty},\\
    \end{EQA}
    with $ \delta_{n,p} $ and $ \delta_{n,p}^{\Wsht} $ from (\ref{Def: Contraction}) for $\Pi$ and $ \Pi^{\Wsht} $, respectively.
\end{theorem}
\begin{remark}
    The error $\Diamond_{NL}$ is due to nonlinearity of the functional and disappears for linear functionals; the error $\Diamond_{GA}$ stands for the Gaussian approximation.
\end{remark}
\begin{remark}
    To make the bound more transparent, consider the Gaussian data, freeze $\St$ and $\phi$, and focus on the dependence on $r$, $p$ and $n$ only. Then
    \begin{EQA}
        \Diamond_{\phi} &\asymp&  \sqrt{\frac{r^3 + r \log^2(n)}{n}}
    \end{EQA}
    for linear functionals, and
    \begin{EQA}
        \Diamond_{\phi} &\asymp&  \sqrt{\frac{p^2 + r^3 + r \log^2(n)}{n}}
    \end{EQA} 
    for nonlinear functionals.
\end{remark}
After that, we just exploit our main result Theorem~\ref{Theorem: posterior independence} to derive the extended version of Theorem~\ref{Theorem: fBvM Wishart} for an arbitrary prior.
\begin{theorem} \label{Theorem: fBvM General}
    Assume the distribution of the data $ \data = (X_1, \ldots, X_{n}) $ fulfills the sample covariance concentration properties (\ref{Eq: hat_delta_n}) and (\ref{ConcentrationCondition2}).
    Suppose the functional $\phi$ satisfies (\ref{Linear expansion}), (\ref{Bound res}).
    Let $ \zeta \sim \ND(0,1) $.
    Then for any prior $ \Pi $ with probability $ 1 - \frac{1}{n} $
    \begin{EQA}
        \sup\limits_{x\in\R} \left| \Ppostbig{\frac{\sqrt{n} \; \left( \phi(\Su) - \phi(\Se) \right)}{\sqrt{2} \| \St^{1/2} \Phib \St^{1/2} \|_2} \leq x} - \Pro(\zeta  \leq x)\right|
        & \lesssim & 
        \Diamond_{\phi} + \Diamond^{\ast}_{\Pi},
    \end{EQA}
    where $ \Diamond_{\phi} $ is defined in (\ref{Eq: Diamond_1}) and $ \Diamond^{\ast}_{\Pi} $ from (\ref{Def: Diamond_ast}) depends on  $ \Pi $.
\end{theorem}
Together with classical results of asymptotic normality of functionals of covariance matrix, these theorems provide a procedure for construction of Bayesian credible sets with frequentist coverage guarantees for the true value of functional at $\St$.

Next, we compare our theorem to the asymptotic functional BvM presented in \cite{Zhou_BvMTFFOCM}.
In that work, a prior distribution is imposed on the precision matrix $ \Ou = \Su^{-1} $. Likewise,  we denote $ \Ot = \St^{-1},\; \Oe = \Se^{-1}$ and so on.
Assume that there exists a set $ A_{n} \subset \SPD$ and a small value $ \delta_{n} = o(1) $ such that
\begin{equation}
    \begin{aligned}
        A_{n} 
    \subset 
    \left\{ \Op \in \SPD: \| \Sp - \St\|_{\infty} \leq \delta_{n} \| \St\|_{\infty} \right\} \, .
    \nonumber
    \end{aligned}
\end{equation}
The condition on approximate linearity of the functional looks like:
\begin{equation}
    \begin{aligned}
     \sup\limits_{\Op \in A_n} \sqrt{n} \| \St^{1/2} \Phib \St^{1/2} \|_{2}^{-1} \left|\phi(\Sp) - \phi(\St) - \Tr\left[ (\Sp-\Se) \Phib \right]\right| = o_{\Pro}(1).
    \label{Linearity}
    \end{aligned}
\end{equation}
The following result takes place.
\begin{theorem}[\cite{Zhou_BvMTFFOCM}, Theorem 2.1] \label{Theorem: BvM linear}
    Under the assumptions of (\ref{Linearity}) and \\
    $ \| \St \|_{\infty} \lor \| \Ot \|_{\infty} = O(1) $, if for a given prior $ \Pi $ the following two conditions are satisfied:
    \begin{enumerate}
        \item $ \Ppost{A_n} = 1 - o_{\Pro}(1) $,
        \item For any fixed $ t \in \R $ holds
                \begin{equation}
                    \begin{aligned}
                        \frac{\int_{A_n} \exp{\left( l_n(\Ou_t)\right)} \; d\Pi(\Ou)}{\int_{A_n} \exp{\left( l_n(\Ou)\right)} \; d\Pi(\Ou)} = 1 + o_{\Pro}(1)
                    \nonumber
                    \end{aligned}
                \end{equation}
                for the perturbed precision matrix
                \begin{equation}
                    \begin{aligned}
                        \Ou_t = \Ou + \frac{\sqrt{2} t}{\sqrt{n} \| \St^{1/2} \Phib \St^{1/2} \|_2}  \Phib,
                    \nonumber
                    \end{aligned}
                \end{equation}
    \end{enumerate} 
    then
        \begin{equation}
            \begin{aligned}
                \sup\limits_{x\in\R} \left| \Ppostbig{\frac{\sqrt{n} \left( \phi(\Su) - \phi(\Se) \right)}{\sqrt{2} \| \St^{1/2} \Phib \St^{1/2} \|_2} \leq x} - \Pro(\zeta \leq x)\right| = o_{\Pro}(1), 
                \nonumber
            \end{aligned}
        \end{equation}
    where $ \zeta \sim \ND(0,1) $.
\end{theorem}
\noindent 
The proof of this result is based on the Laplace transform approach, which by no means allows to recover the rate of approximation error in Kolmogorov distance. Moreover, the proof relies heavily on the ``shift absorption'' condition which is not straightforward, and the authors consider only the case of the Gaussian data $ X_1, \ldots, X_n $.
\begin{remark}
    \cite{Zhou_BvMTFFOCM} also derives similar results for functionals of precision matrix $\Ou$. Our technique allows to do that as well. However, the theory for this case is even simpler, because the conjugate prior to the multivariate Gaussian distribution parametrized by the precision matrix is the Wishart distribution, which is more straightforward to work with. Thus, we focus on covariance matrices.
\end{remark}

\subsubsection{``Rank-adjusted'' non-conjugate priors and dimension-free result} \label{Rank-adjusted Prior}
Note that the dimension-free rates for the conjugate prior obtained in Theorem \ref{Theorem: fBvM Wishart} suffer when we pass to a non-conjugate prior in Theorem \ref{Theorem: fBvM General} as the dimension of the space $p$ appears in $\Diamond_{\Pi}^*$. However, in order to get a dimension-free bound for linear functionals, we can make use of potential low rank of $\Phib$. In particular, we introduce ``rank-adjusted'' non-conjugate priors, which in practice allow to reduce the computational complexity and improve the quality of the procedure allowing the dimension $p$ to be arbitrarily large.

Once we have $\phi(\Su) = \Tr[\Phib \Su]$ with $\Phib$ of rank $r$, we compute its spectral decomposition:
\begin{equation}
    \begin{aligned}
        \Phib = \U\Psib \U^{\top},
        \nonumber
    \end{aligned}
\end{equation}
where $\Psib \in \R^{r\times r}$ is diagonal and $\U\in\R^{p\times r}$ satisfies $\U^{\top} \U = \I_r$. As previously, we impose a prior distribution $\widetilde{\Pi}$ on covariance matrices, but now on the set $\mathbbm{S}^{r}_{+}$ rather than $\SPD$:
\begin{equation}
    \begin{aligned}
        \Xiu \sim \widetilde{\Pi}.
        \nonumber
    \end{aligned}
\end{equation}
Its counterpart in $\SPD$ 
\begin{equation}
    \begin{aligned}
        \Su \eqdef \U \Xiu \U^{\top}
        \nonumber
    \end{aligned}
\end{equation}
is distributed according to some prior $\Pi$ on $\SPD$ which we call ``rank-adjusted'' prior.
After the data $X_1, \ldots, X_n \in \R^p$ are observed, we modify them: for $j = 1, \ldots, n$
\begin{equation}
    \begin{aligned}
        Y_j \eqdef \U^{\top} X_j \in \R^{r},
        \nonumber
    \end{aligned}
\end{equation}
which clearly has the true covariance $\Xit \eqdef \U^{\top} \St \U$ and implies the sample covariance $\Xie \eqdef \U^{\top} \Se \U$. 
It's trivial to see
\begin{equation}
    \begin{aligned}
        &\Tr[\Phib \St] = \Tr[\Psib \Xit],\;\;\; \Tr[\Phib \Se] = \Tr[\Psib \Xie],\;\;\; \Tr[\Phib \Su] = \Tr[\Psib \Xiu],\\
        &\| \St^{1/2} \Phib \St^{1/2} \|_2 = \| {\Xit}^{1/2} \Psib {\Xit}^{1/2}\|_2,
        \nonumber
    \end{aligned}
\end{equation}
so we have the following equality in distribution
\begin{equation}
    \begin{aligned}
        \left(\sqrt{n}\, \frac{\Tr[\Phib \Su] - \Tr[\Phib \Se]}{\sqrt{2} \| \St^{1/2} \Phib \St^{1/2}\|_2} \;\Bigg|\; \data\right)
        \;\;\stackrel{d}{=}\;\;
        \left(\sqrt{n}\, \frac{\Tr[\Psib \Xiu] - \Tr[\Psib \Xie]}{\sqrt{2} \| {\Xit}^{1/2} \Psib {\Xit}^{1/2}\|_2} \;\Bigg|\; \boldsymbol{Y}^n\right),
        \nonumber
    \end{aligned}
\end{equation}
or
\begin{equation}
    \begin{aligned}
    \Pi\left(\sqrt{n}\, \frac{\Tr[\Phib \Su] - \Tr[\Phib \Se]}{\sqrt{2} \| \St^{1/2} \Phib \St^{1/2}\|_2} \leq x\;\Bigg|\; \data\right)
        \,=\,
        \widetilde{\Pi}\left(\sqrt{n}\, \frac{\Tr[\Psib \Xiu] - \Tr[\Psib \Xie]}{\sqrt{2} \| {\Xit}^{1/2} \Psib {\Xit}^{1/2}\|_2} \leq x\;\Bigg|\; \boldsymbol{Y}^n\right)
        \nonumber
    \end{aligned}
\end{equation}
for all $x\in\R$ with probability $1$.
Hence, Theorem \ref{Theorem: fBvM General} applied to the data $Y_1, \ldots, Y_n$ with true covariance $\Xit$ and sample covariance $\Xie$, prior $\widetilde{\Pi}$ on $\Xiu$ and linear functional described by $\Psib$ yields the approximation by the standard Gaussian distribution with the error in terms of the rank $r$ rather than the dimension $p$.

		\subsection{Spectral projectors of covariance matrix} \label{Subsection: projectors}
			In this Subsection we follow \cite{Silin}. Define the empirical projector on the $s$-th eigenspace of dimension $\ms_{s}$ that comes from the sample covariance $\Se$ and estimates the true projector $\Ptr$:
\begin{EQA}
    \Per 
    & = &
    \sum\limits_{j\in\Deltas_{s}} \uh_j {\uh_j}^{\T}.
\end{EQA}
Similarly, for the covariance matrix $\Su$, which is generated from the pseudo-posterior distribution, we write
\begin{EQA}
    \Pur
    & = &
    \sum\limits_{j\in\Deltas_{s}} u_j {u_j}^{\T},
\end{EQA}
where $u_1, \ldots, u_p$ are the eigenvectors of $\Su$ associated with eigenvalues ordered in descending order.

More generally, we can pick any of the $q$ eigenspaces by fixing the set 
$\subspaceset \subset \{ 1, \ldots, q \}$. To avoid technical difficulties, we restrict $\subspaceset$ to be of the form
\begin{EQA}[c]
    \subspaceset = \{ s^{-}, s^{-}+1, \ldots, s^{+}\}
\end{EQA}
for some $ 1 \leq s^{-} \leq s^{+} \leq q$.
Below it is useful to introduce the subset of indices 
\begin{EQA}
    \IndS 
    & \eqdef &
    \bigl\{ k \colon k \in \Deltas_s, s \in \subspaceset \bigr\} \,.
\end{EQA}
Now we define projector onto the direct sum of the eigenspaces associated with  $\Ptr$ for all $s \in \subspaceset$:
\begin{equation}
    \begin{aligned}
        \Pt \eqdef \sum\limits_{s \in \subspaceset} \Ptr = \sum\limits_{k \in \IndS} \us_k {\us_k}^{\T}.
    \nonumber
    \end{aligned}
\end{equation}
Its empirical counterpart is given by
\begin{equation}
    \begin{aligned}
        \Pe \eqdef \sum\limits_{s \in \subspaceset} \Per = \sum\limits_{k \in \IndS} \uh_k \uh_k^{\T},
    \nonumber
    \end{aligned}
\end{equation}
while the Bayesian version derived from pseudo-posterior reads as
\begin{equation}
    \begin{aligned}
        \Pu \eqdef \sum\limits_{s \in \subspaceset} \Pur = \sum\limits_{k \in \IndS} u_k u_k^{\T}.
    \nonumber
    \end{aligned}
\end{equation}
Such objects are of greater interest, since, for instance, when $\subspaceset = \{ 1, \ldots, d\}$ for some $d < q$, then $\Pe$ is exactly what is recovered by PCA. It turns out, that the random quantity from Bayesian world
$n\| \Pu - \Pe\|_2^2$ mimics the distribution of $n\| \Pe - \Pt\|_2^2$ (see again \cite{Silin}) and therefore is of interest in PCA-related problems.

The total rank of $\Pt$ is $\ms_{\subspaceset} = \sum\limits_{s \in \subspaceset} \ms_s$. 
Further, we need to introduce spectral gap for the set $\subspaceset$:
\begin{EQA}[c]
        \gs_{\subspaceset} \eqdef \begin{cases}
        \mus_{s^+} - \mus_{s^+ +1}, & \text{ if } s^- = 1;\\
        \mus_{s^- - 1} - \mus_{s^-}, & \text{ if } s^+ = q; \\
        \left( \mus_{s^--1}-  \mus_{s^-}\right) \;\land\; \left( \mus_{s^+} - \mus_{s^+ +1} \right), & \text{ otherwise}. 
        \end{cases}
\end{EQA}
To describe the posterior distribution of $n\| \Pu - \Pe \|_{\Fr}^2$, 
we introduce the following block-diagonal matrix $\Gammas_{\subspaceset}$ of size 
$\ms_{\subspaceset} (p - \ms_{\subspaceset}) \times \ms_{\subspaceset} (p - \ms_{\subspaceset})$:
\begin{equation}
    \begin{aligned}
        \Gammas_{\subspaceset} & \eqdef  \diag\left(\Gamma_{\subspaceset}^s\right)_{s \in \subspaceset}, \\
        \Gamma_{\subspaceset}^s & \eqdef  \diag\left(\Gamma^{s,k}\right)_{k \notin \subspaceset},\\
        \Gamma^{s,k} & \eqdef  \frac{2\mus_{s} \mus_{k}}{(\mus_{s} - \mus_{k})^2} \cdot \I_{\ms_s \ms_k}, \;\;\;\; s \in \subspaceset, k \notin \subspaceset.
        \label{Formula: Gamma_r}
    \end{aligned}
\end{equation}

Before formulating the result for an arbitrary prior, we consider the Inverse Wishart prior.
\begin{theorem}[\cite{Silin}, Theorem 2.1] \label{Theorem: BvM projector conjugate}
    Assume the distribution of the data $\data = (X_{1}, \ldots, X_{n})$ fulfills the sample covariance concentration property (\ref{Eq: hat_delta_n}) with $ \hdelta_{n,p} $ satisfying
    \begin{equation}
            \begin{aligned}
            & \hdelta_{n,p} \leq \frac{\gs_{\subspaceset}}{4\| \St \|_{\infty}} \;\land\; \frac{\reff(\St)}{p}.
            \nonumber
            \end{aligned}
        \end{equation}
Consider the prior $\Pi^{\mathcal{W}}$ given by the Inverse Wishart distribution
$\IWsht_{p}(\G,\, p+b-1)$. Let $\xi \sim \ND(0, \Gammas_{\subspaceset})$ with 
$\Gammas_{\subspaceset}$ defined by (\ref{Formula: Gamma_r}). 
Then with probability $1 - 1/n$
\begin{EQA}[c]
    \sup_{x\in\R }\left| 
    	\PpostW{ n\| \Pu - \Pe \|_{2}^{2} \leq x } - \Pro( \| \xi \|^{2} \leq x)  
    \right| 
    \lesssim  
    \Diamond_P,
\end{EQA}
    where
\begin{equation}
    \begin{aligned}
    \Diamond_P = \Diamond_P(n,p,\St)
    \eqdef 
    \frac{\Diamond_{1} + \Diamond_{2} + \Diamond_{3}}
         {\| \Gammas_{\subspaceset} \|_2^{1/2} \left( \| \Gammas_{\subspaceset} \|_2^2 - \| \Gammas_{\subspaceset} \|_{\infty}^2 \right)^{1/4}} 
             + \frac{1}{n}\, .
    \label{Eq: Diamond}
    \end{aligned}
\end{equation}
    The terms $ \Diamond_{1} $ through $ \Diamond_{3} $ can be described as
\begin{EQA}
    \Diamond_{1} 
    & \eqdef & 
    \left\{ 
        (\log(n) + p) \left( \,\left( 1 + \frac{l^*_{\subspaceset}}{\gs_{\subspaceset}} \right) \frac{\sqrt{\ms_{\subspaceset}} \| \St \|_{\infty}}{\gs_{\subspaceset}} + \ms_{\subspaceset} \right) \| \St \|_{\infty} + \ms_{\subspaceset} \,\| \G \|_{\infty}
    \right\}  \\
    && \hspace{6.95cm} \times \frac{\ms_{\subspaceset}\,\|\St\|_{\infty}}{{\gs_{\subspaceset}}^{2}}  \sqrt{\frac{\log(n) + p}{n}}
    ,\\
    \Diamond_{2}
    & \eqdef & 
    \frac{ \| \St \|_{\infty} \,\left(\ms_{\subspaceset} \, \| \St \|_{\infty}^{2} \land \Tr\left({\St}^{2}\right) \right) }{{\gs_{\subspaceset}}^3} \, p\left( \hdelta_{n,p} + \frac{p}{n} \right),\\
    \Diamond_{3} 
    & \eqdef &
    \frac{ {(\ms_{\subspaceset})}^{3/2} \| \St \|_{\infty} \, \Tr(\St)}{{\gs_{\subspaceset}}^{2}}
    \sqrt{\frac{\log(n)}{n}}.
\end{EQA}
\end{theorem}
\begin{remark}[\cite{Silin}, Remark 2.1]
The bound (\ref{Eq: Diamond}) can be made more transparent 
if we fix $ \St $ and focus on the dependence on $ p, n, \hdelta_{n,p} $ and the desired subspace dimension $ \ms_{\subspaceset} $ only (freezing the eigenvalues, the spectral gaps and multiplicities of the eigenvalues):
        \begin{equation}
            \begin{aligned}
                \Diamond_P \;\asymp\; \sqrt{\frac{(\ms_{\subspaceset})^3\,(p^3 + \log^3(n))}{n}} \,+\, \ms_{\subspaceset}\,p\,\hdelta_{n,p},
            \nonumber
            \end{aligned}
        \end{equation}
or, in the sub-Gaussian case,
    \begin{equation}
            \begin{aligned}
                \Diamond_P \;\asymp\; \sqrt{\frac{(\ms_{\subspaceset})^3\,(p^3 + \log^3(n))}{n}}\,.
            \nonumber
            \end{aligned}
        \end{equation}
Moreover, in the case of spiked covariance model we expect $\| \Gammas_{\subspaceset} \|_2$ to behave as $\sqrt{p\,\ms_{\subspaceset}}$ which improves the previous bounds to
\begin{equation}
            \begin{aligned}
                \Diamond_P \;\asymp\; \sqrt{\frac{(\ms_{\subspaceset})^2\,(p^2 + \log^3(n)/p)}{n}} \,+\, \sqrt{\ms_{\subspaceset}\,p}\,\hdelta_{n,p},
            \nonumber
            \end{aligned}
        \end{equation}
and
    \begin{equation}
            \begin{aligned}
                \Diamond_P \;\asymp\; \sqrt{\frac{(\ms_{\subspaceset})^2\,(p^2 + \log^3(n)/p)}{n}},
            \nonumber
            \end{aligned}
        \end{equation}
respectively.
\end{remark}

From Theorem \ref{Theorem: BvM projector conjugate} and Theorem \ref{Theorem: posterior independence} we derive the following.
\begin{theorem} \label{Theorem: BvM projector finite}
    Assume the distribution of the data $ \data = (X_1, \ldots, X_{n}) $ fulfills the sample covariance concentration property (\ref{Eq: hat_delta_n}) with $ \hdelta_{n,p} $ satisfying
    \begin{equation}
            \begin{aligned}
            & \hdelta_{n,p} \leq \frac{\gs_{\subspaceset}}{4\| \St \|_{\infty}} \;\land\; \frac{\reff(\St)}{p}.
            \nonumber
            \end{aligned}
        \end{equation}
    Let $ \xi \sim \ND(0, \Gammas_{\subspaceset}) $ with $ \Gammas_{\subspaceset} $ defined by (\ref{Formula: Gamma_r}).
    Then for any prior $ \Pi $ with probability $ 1 - 1/n $
    \begin{EQA}
            \sup\limits_{x\in\R }\left| \Ppost{ n\| \Pu - \Pe \|_2^2 \leq x } - \Pro( \| \xi \|^2 \leq x)  \right| 
            & \lesssim & 
            \Diamond_P + \Diamond^{\ast}_{\Pi},
    \end{EQA}
    where $ \Diamond_P $ is defined in (\ref{Eq: Diamond}) and $ \Diamond^{\ast}_{\Pi} $ from (\ref{Def: Diamond_ast}) depends on $ \Pi $.
\end{theorem}
As for the functionals of covariance matrix, this result can be applied for building of sharp confidence sets for the true projector $ \Pt $.
See again \cite{Silin}, Corollary 2.3 and the paragraph before for the detailed description of the procedure.
    \section{Discussion and conclusions}
    \label{Section: Discussion}
        Here we discuss the approximation error rates obtained in previous chapters and compare with previous results from the literature.
We consider linear functionals such as
\begin{itemize}
    \item entry ${\Su}_{ij}$ with $\Phib = \frac{e_i e_j^{\top} + e_j e_i^{\top}}{2}$ ($e_i$ denotes $i$-th vector of the standard basis in $\R^p$),
    \item quadratic form $v^{\top} \Su v$ with $\Phib = vv^{\top}$,
    \item trace $\Tr[\Su]$ with $\Phib=\I_p$,
\end{itemize}
as well as nonlinear functionals: eigenvalues, log-determinant, and finally, spectral projectors.

We fix such characteristics of the underlying true covariance matrix such as $\| \St \|_{\infty}$, $\| \St^{-1} \|_{\infty}$, spectral gaps of $\St$ and study the dependence on the sample size $n$, dimension of the data $p$, the rank of the functional $r$ and some additional parameters in a couple of cases, dropping logarithmic terms. For the purposes of comparison, we consider the case of the Gaussian data. 

In the following table, the column ``conjugate'' corresponds to the rates obtained for the conjugate prior in Theorem~\ref{Theorem: fBvM Wishart} and Theorem~\ref{Theorem: BvM projector conjugate}, the column ``non-conjugate'' corresponds to the rates obtained for an arbitrary non-conjugate prior in Theorem~\ref{Theorem: fBvM General} and Theorem~\ref{Theorem: BvM projector finite}, and the column ``non-conj. rank-adj.'' corresponds to the class of rank-adjusted non-conjugate priors defined in Subsection~\ref{Rank-adjusted Prior}.
\\
\begin{center}
\begin{tabular}{|l||c|c|c|c|}
\hline
  &  conjugate & non-conjugate & non-conj. \\
  & & & rank-adj. \\ \hline \hline
 Linear functionals & $\sqrt{r^3/n}$ & $\sqrt{p^5/n}$ & $\sqrt{r^5/n}$   \\ \hline
 - Entry ${\Su}_{ij}$ & $\sqrt{1/n}$ & $\sqrt{p^5/n}$ & $\sqrt{1/n}$ \\ \hline
 - Quadratic form $v^{\top} \Su v$ & $\sqrt{1/n}$  & $\sqrt{p^5/n}$ & $\sqrt{1/n}$ \\ \hline
 - Trace $\Tr[\Su]$ & $\sqrt{p^3/n}$ & $\sqrt{p^5/n}$ & $\sqrt{p^5/n}$\\ \hline \hline
 Nonlinear functionals &  $\sqrt{(p^2+r^3)/n}$  & $\sqrt{p^5/n}$ & $\times$ \\ \hline
 - Eigenvalues (of multiplicity $m$) & $\sqrt{(p^2 + m^3)/n}$ & $\sqrt{p^5/n}$ & $\times$ \\ \hline
 - Log-determinant & $\sqrt{p^3/n}$ & $\sqrt{p^5/n}$ & $\times$ \\ \hline \hline
 Spectral projectors (of rank $m$) & $\sqrt{(p^3 m^3)/n}$ & $\sqrt{(p^5+p^3 m^3)/n}$ & $\times$ \\ \hline
\end{tabular}
\end{center}
\vspace{0.4cm}
One may compare this table to the one provided in \cite{Zhou_BvMTFFOCM}, Section 5.2. Some of the results for the conjugate prior are significantly improved, and now for all of the presented functionals the sharp rates from frequentist asymptotic normality results (see \cite{Zhou_BvMTFFOCM}, Section 5.1) match the rates for the conjugate prior. The rates for non-conjugate priors require $p^5 \ll n$; we already discussed in Remark~\ref{Remark BvM} that this is the payment for simplicity of the proof technique and generality of prior and underlying data distribution. Another explanation is that Theorem~\ref{Theorem: posterior independence} is too general in a sense that it provides bound of the total variation distance between posterior distributions of $\Su$ without making use of specific application, such as functionals or spectral projectors. We see that in case of spectral projectors this payment for the non-conjugate prior may be even smaller than the error of approximation for the conjugate prior, if the rank of the projector $m \gtrsim p^{2/3}$. Observe that, compared to the general non-conjugate priors, the ``rank-adjusted'' non-conjugate priors described in Subsection~\ref{Rank-adjusted Prior} provide significantly better rates for low-rank linear functionals.

Note once again that \cite{Zhou_BvMTFFOCM} is able only to specify regimes in which the convergence of posterior to the Gaussian distribution occurs, under the assumption of Gaussian data and the specific condition on a prior. In contrast, though our rates in case of non-conjugate priors turn out to be worse in terms of the required relation between $p$ and $n$, they are explicit and work for finite sample and much wider classes of data distributions under simpler prior assumption.

	 \section{Main proofs} \label{Section: Proofs}
	 	\subsection{Proof of Theorem \ref{Theorem: posterior independence}}
	 		\paragraph{Step 1 } \textit{``Localization''}.\\
Fix arbitrary prior $\Pi$ (in particular, we may take $\Pi^{\Wsht}$).
Posterior measure of a set $\A \subset \SPD$ is
\begin{EQA}
    \Pi(\A | \data) 
    & = & 
    \frac{\int_{\A} \exp{(l_n(\Su))} \;\Pi(\Su) \;d\Su }{\int_{\SPD} \exp{(l_n(\Su))} \;\Pi(\Su) \;d\Su} = \frac{\tau(\A)}{\tau(\SPD)},
\end{EQA}
where for shortness we introduce
\begin{EQA}
    \tau(\A) 
    & \eqdef & 
    \int_{\A} \exp{(l_n(\Su))} \; \Pi(\Su) \;d\Su.
\end{EQA}
Observe that due to (\ref{Def: Contraction}), since $ \overline{\delta}_{n,p} \geq \delta_{n,p} $, we have
\begin{equation}
    \begin{aligned}
	\frac{\tau(\vic(\overline{\delta}_{n,p}))}{\tau(\SPD)} 
	\geq 
	 1 - \frac{1}{n}
	 \label{tau conc}
	 \end{aligned}
\end{equation}
with probability $1 - 1/n$.
Consider \textit{``localized''} posterior defined by
\begin{EQA}
    \Pi_{loc} (\A | \data) 
    & \eqdef &
    \frac{\int_{\A \cap \vic(\overline{\delta}_{n,p})} \exp{(l_n(\Su))} \;\Pi(\Su) \;d\Su }{\int_{\vic(\overline{\delta}_{n,p})} \exp{(l_n(\Su))} \;\Pi(\Su) \;d\Su} = 
    \frac{\tau(\A \cap \vic(\overline{\delta}_{n,p}))}{\tau(\vic(\overline{\delta}_{n,p}))} \,. 
\end{EQA}
It straightforwardly follows from (\ref{tau conc}) that
\begin{equation}
    \begin{aligned}
    \sup\limits_{\A \subset \SPD} \left| \Pi_{loc}(\A \cond \data) - \Pi(\A \cond \data) \right| 
     \leq 
    \frac{2}{n}.
    \label{Bound: localization}
    \end{aligned}
\end{equation}
with probability $1 - 1/n$.
In particular, Lemma~\ref{Contraction} and the fact that $ \overline{\delta}_{n,p} \geq \delta_{n,p}^{\Wsht} $ imply similar bound for the Inverse Wishart prior:
\begin{equation}
    \begin{aligned}
    \sup\limits_{\A \subset \SPD} \left| \Pi^{\Wsht}_{loc}(\A \cond \data) - \Pi^{\Wsht}(\A \cond \data) \right| 
    \leq 
    \frac{2}{n}.
    \label{Bound: localization Wishart}
\end{aligned}
\end{equation}

\paragraph{Step 2 } \textit{``Flatness''}.\\
Consider uniform prior $\Pi^U$ over $\vic(\overline{\delta}_{n,p})$, which will play role of a bridge between the priors $\Pi$ and $\Pi^{\Wsht}$. 
The corresponding posterior is given by
\begin{EQA}
    \Pi^U (\A \cond \data) 
    & = &
    \frac{\int_{\A \cap \vic(\overline{\delta}_{n,p}) } \exp{(l_n(\Su))} \;d\Su }{\int_{\vic(\overline{\delta}_{n,p})} \exp{(l_n(\Su))} \;d\Su}.
\end{EQA}
Recalling the definition of flatness (\ref{Def: flatness}), it is easy to show that
\begin{equation}
    \begin{aligned}
   	\sup\limits_{\A \subset \SPD} \left| \Pi_{loc}(\A \cond \data) - \Pi^U(\A \cond \data) \right| 
    \leq 
    4\, \fl(\overline{\delta}_{n,p})
    \label{Bound: uniform}
\end{aligned}
\end{equation}
with probability $1$. 
The same applies to the Inverse Wishart prior:
\begin{equation}
    \begin{aligned}
   	\sup\limits_{\A \subset \SPD} \left| \Pi^{\Wsht}_{loc}(\A \cond \data) - \Pi^U(\A \cond \data) \right| 
    \leq 
    4\, \fl^{\Wsht}(\overline{\delta}_{n,p})
    \label{Bound: uniform Wishart}
    \end{aligned}
\end{equation}
with probability $1$. 

Applying the triangle inequality to (\ref{Bound: localization}), (\ref{Bound: localization Wishart}), (\ref{Bound: uniform}) and (\ref{Bound: uniform Wishart}), we derive the desired result.
	 	\subsection{Proof of Theorem \ref{Theorem: fBvM Wishart}}
			Since $\Su \sim \Wsht_p^{-1}(\G,\,p + b -1)$, then from the prior density $\Pi^{\Wsht}$ and the Gaussian log-likelihood $l_n$ one may compute the posterior density and notice that
\begin{equation}
    \begin{aligned}
        \Su \,|\, \data \sim \Wsht_p^{-1}(n\Se + \G,\, n+p+b-1).
        \nonumber
    \end{aligned}
\end{equation}
Throughout the proof we will be using the properties of the Inverse Wishart distribution presented in Lemma~\ref{Lemma: Properties} from Appendix~\ref{Appendix: A}.

Let us get rid of possible non-linearity of our functional $\phi(\cdot)$.
Write the representation (\ref{Linear expansion}) in two ways:
\begin{EQA}
	 \phi(\Su) - \phi(\St) 
	 & = &
	 \Tr[\Phib\,(\Su-\St)] + \res(\Su,\St),\\
	 \phi(\Se) - \phi(\St) 
	 & = &
	 \Tr[\Phib\,(\Se-\St)] + \res(\Se,\St).
\end{EQA}
Thus, subtracting the latter equality from the former,
\begin{EQA}
	 \phi(\Su) - \phi(\Se) 
	 & = &
	 \Tr[\Phib(\Su-\Se)] + \res(\Su,\St) - \res(\Se,\St).
\end{EQA}
Define
\begin{equation}
\begin{aligned}
	\label{Delta1}
    \Delta_1 \eqdef  \sqrt{n}\, {\delta_{n,p}^{\Wsht}}^2 \,\CONST_{\phi}(\St)\, \| \St \|_{\infty}^2,
	\end{aligned}
\end{equation}
and
\begin{equation}
\begin{aligned}
    \Delta_2 \eqdef  \sqrt{n}\, \hdelta_{n,p}^2 \,\CONST_{\phi}(\St)\, \| \St \|_{\infty}^2
    \label{Delta2}
	\end{aligned}
\end{equation}
with $ \delta_{n,p} $ and $ \delta_{n,p}^{\Wsht} $ from (\ref{Def: Contraction}) for $\Pi$ and $ \Pi^{\Wsht} $, respectively.
Then the assumption (\ref{Bound res}) yields
\begin{equation}
\begin{aligned}
	\label{Prob_Delta1}
	\PpostW{\sqrt{n}\, |\res(\Su,\St)| > \Delta_1} \leq \frac{1}{n}\;\;\; \text{ with probability }1-\frac{1}{n},
\end{aligned}
\end{equation}
and
\begin{equation}
\begin{aligned}
	\Pro(\sqrt{n}\, |\res(\Se,\St)| > \Delta_2) \leq \frac{1}{n}.
	\label{Prob_Delta2}
\end{aligned}
\end{equation}
Further, we elaborate on the linear part $\sqrt{n}\,\Tr \left[ \Phib\,(\Su - \Se)\right]$. Recall that
\begin{equation}
    \begin{aligned}
        \St^{1/2} \Phib \St^{1/2} = \V \D \V^{\top} \;\text{with}\; \V\in\R^{p\times r}, \V^{\top} \V = \I_r \text{ and }\D\in\R^{r\times r} \text{ diagonal}.
        \nonumber
    \end{aligned}
\end{equation}
Introduce for shortness
\begin{equation}
    \begin{aligned}
        &n_{r,b} \eqdef n+r+b-1,\;\;\;\;
        n_{corr} \eqdef \frac{n_{r,b}^{3/2}}{\sqrt{n}},\\
        &\Aa \eqdef \V^{\top} \St^{-1/2} \left(\frac{n\Se + \G}{n_{corr}}\right) \St^{-1/2} \V \,\in \mathbbm{S}^{r}_{+},\\
        &\Bb \eqdef \Aa^{-1/2}\V^{\top} \St^{-1/2} \Su \St^{-1/2} \V \Aa^{-1/2} \,\in \mathbbm{S}^{r}_{+},\\
        &\Gg \eqdef \V^{\top} \St^{-1/2} \G \St^{-1/2} \V.
        \nonumber
    \end{aligned}
\end{equation}
Note that by Lemma \ref{Lemma: Properties} (ii) it holds
\begin{equation}
    \begin{aligned}
        \frac{\Bb}{n_{corr}} \sim \Wsht_r^{-1}(\I_r,\, n_{r,b}).
        \nonumber
    \end{aligned}
\end{equation}
Also we need to show that $\Aa$ is close to the identity matrix. The following auxiliary lemma quantifies this statement. We will use this result later in the proof.
\begin{lemma} \label{Lemma: A}
It holds
\begin{equation}
    \begin{aligned}
        \| \Aa^{1/2} - \I_r \|_{\infty} \leq \| \Aa - \I_r \|_{\infty} \leq \wdelta_{n,r} + \frac{r}{n_{r,b}} + \frac{\| \Gg \|_{\infty}}{n_{r,b}}
        \nonumber
    \end{aligned}
\end{equation}
with probability $1-1/n$.
Consequently,
\begin{equation}
    \begin{aligned}
        \| \Aa \|_{\infty} \lesssim 1 + \frac{\| \Gg \|_{\infty}}{n_{r,b}}
        \nonumber
    \end{aligned}
\end{equation}
with probability $1-1/n$.
\end{lemma}
\begin{proof}
    Let us prove the first inequality. Let $\lambda$ be an eigenvalue of $\Aa - \I_r$, and since $\Aa \in \mathbbm{S}^{r}_{+}$ then $\lambda \geq -1$. The corresponding eigenvalue of $\Aa^{1/2}-\I_r$ is $(\sqrt{1+\lambda}-1)$. Trivial inequality $|\sqrt{1+\lambda}-1| \leq |\lambda|$ concludes the first part of the proof.
    
    Now we prove the second inequality. We have
    \begin{equation}
        \begin{aligned}
            \| \Aa - \I_r \|_{\infty} 
            &= \left\| \V^{\top} \St^{-1/2} \left(\frac{n\Se + \G}{n_{corr}}\right) \St^{-1/2} \V - \I_r \right\|_{\infty}\\
            &\leq \| \V^{\top} \St^{-1/2} \left(\frac{n\Se - n_{corr}\St }{n_{corr}}\right) \St^{-1/2} \V \|_{\infty} + \frac{\| \Gg \|_{\infty}}{n_{corr}}\\
            &\leq \| \V^{\top} \St^{-1/2} (\Se - \St) \St^{-1/2} \V \|_{\infty} + \frac{n_{corr} - n}{n_{corr}} + \frac{\| \Gg \|_{\infty}}{n_{r,b}}\\
            &\leq \wdelta_{n,r} + \frac{r+b}{n_{r,b}} + \frac{\| \Gg \|_{\infty}}{n_{r,b}},
            \nonumber
        \end{aligned}
    \end{equation}
    where the last inequality holds with probability $1-1/n$.
    Moreover,
    \begin{equation}
        \begin{aligned}
            \|\Aa\|_{\infty} \leq 1 + \| \Aa - \I_r \|_{\infty},
            \nonumber
        \end{aligned}
    \end{equation}
    and the ``consequently'' part of the claim follows assuming without loss of generality that $\wdelta_{n,r} < 1$ (otherwise the statement of the whole Theorem is trivial).
\end{proof}
Our first approximation result is as follows.
\begin{lemma} \label{Lemma: Third}
    For the remainder $\mathcal{R}_3$ defined as
    \begin{equation}
        \begin{aligned}
            \mathcal{R}_3 \eqdef \sqrt{n}\Tr\left[ \Phib\,(\Su-\Se)\right] - \sqrt{n} \Tr\left[\Aa^{1/2}\D\Aa^{1/2} \left(\Bb - \frac{n_{corr}}{n_{r,b}}\,\I_r\right)\right]
            \nonumber
        \end{aligned}
    \end{equation}
    and 
    \begin{equation}
        \begin{aligned}
            \Delta_3 \stackrel{\operatorname{def}}{\asymp} 
            \| \D\|_1 \,\left( \frac{\|\Gg\|_{\infty}}{\sqrt{n_{r,b}}} + \frac{r}{\sqrt{n_{r,b}}}\right)
            \label{Delta3}
        \end{aligned}
    \end{equation}
   holds
   \begin{equation}
        \begin{aligned}
            \Pi^\Wsht\big( |\mathcal{R}_3| \leq \Delta_3 \;|\;\data\big) = 1
            \nonumber
        \end{aligned}
    \end{equation}
    with probability $1-1/n$.
    
\end{lemma}
\begin{proof}
    Using the cyclic property of the trace, the spectral decomposition of $\St^{1/2} \Phib \St^{1/2}$ and the definitions of $\Aa$ and $\Bb$, we have
    \begin{equation}
        \begin{aligned}
            \sqrt{n}\Tr\left[ \Phib\,(\Su-\Se)\right] &=
            \sqrt{n}\Tr\left[ \St^{1/2} \Phib \St^{1/2}\,\St^{-1/2}(\Su-\Se)\St^{-1/2}\right]\\
            &=\sqrt{n}\Tr\left[ \V\D\V^{\top}\cdot\St^{-1/2}(\Su-\Se)\St^{-1/2}\right]\\
            &=\sqrt{n}\Tr\left[ \D\cdot\V^{\top}\St^{-1/2}(\Su-\Se)\St^{-1/2}\V\right]\\
            &=\sqrt{n}\Tr\left[ \D\,(\Aa^{1/2}\Bb\Aa^{1/2} - \V^{\top}\St^{-1/2}\Se\St^{-1/2}\V)\right].
            \nonumber
        \end{aligned}
    \end{equation}
    Now we approximate $\V^{\top}\St^{-1/2}\Se\St^{-1/2}\V$ by $\frac{n_{corr}}{n_{r,b}}\,\Aa$ in spectral norm:
    \begin{equation}
        \begin{aligned}
            &\left\| \V^{\top}\St^{-1/2}\Se\St^{-1/2}\V - \frac{n_{corr}}{n_{r,b}}\,\Aa\right\|_{\infty} =\\
            &\hspace{1cm} = \left\| \V^{\top}\St^{-1/2}\left( \Se - \frac{n\Se+\G}{n_{r,b}}\right)\St^{-1/2}\V \right\|_{\infty}\\
            &\hspace{1cm} \leq \frac{\|\Gg\|_{\infty}}{n_{r,b}} + \frac{r+b}{n_{r,b}}\,\| \V^{\top}\St^{-1/2}\Se\St^{-1/2}\V \|_{\infty}\\
            &\hspace{1cm} \leq \frac{\|\Gg\|_{\infty}}{n_{r,b}} + \frac{r+b}{n_{r,b}}\,(1 + \wdelta_{n,r}),
            \nonumber
        \end{aligned}
    \end{equation}
    where the last inequality holds with probability $1-1/n$.
    Hence,
    \begin{equation}
        \begin{aligned}
            &\left| \sqrt{n}\Tr\left[ \Phib\,(\Su-\Se)\right] - \sqrt{n} \Tr\left[\D\,(\Aa^{1/2} \Bb \Aa^{1/2} - \frac{n_{corr}}{n_{r,b}}\,\Aa)\right] \right| \leq\\
            &\hspace{1cm}\leq\sqrt{n} \,\| \D\|_1 \, \left\| \V^{\top}\St^{-1/2}\Se\St^{-1/2}\V - \frac{n_{corr}}{n_{r,b}}\,\Aa\right\|_{\infty}\\
            &\hspace{1cm}\lesssim \| \D\|_1 \,\left( \frac{\|\Gg\|_{\infty}}{\sqrt{n_{r,b}}} + \frac{r+b}{\sqrt{n_{r,b}}}\right),
            \nonumber
        \end{aligned}
    \end{equation}
    where we again used the assumption that $\wdelta_{n,r} < 1$.
    This concludes the proof of the lemma.
\end{proof}
    Introduce 
    \begin{equation}
        \begin{aligned}
            \Ee = \frac{1}{n_{r,b}} \sum\limits_{j=1}^{n_{r,b}} Z_j Z_j^{\top} - \I_r,
            \nonumber
        \end{aligned}
    \end{equation}
    where $\{ Z_j \}_{j=1}^{n_{r,b}}\,|\,\data \stackrel{i.i.d.}{\sim} \mathcal{N}_r(0, \I_r)$. From Lemma \ref{Lemma: Properties} (iii) we know that
    \begin{equation}
        \begin{aligned}
            \frac{n_{r,b}}{n_{corr}}\,\Bb\;\Big|\;\data \eqdist (\I_r + \Ee)^{-1}.
            \nonumber
        \end{aligned}
    \end{equation}
    and therefore
    \begin{equation}
        \begin{aligned}
            \Bb -\frac{n_{corr}}{n_{r,b}}\,\I_r \;\Big|\;\data \eqdist \frac{n_{corr}}{n_{r,b}}\, \left[(\I_r + \Ee)^{-1} - \I_r\right].
            \label{Distribution}
        \end{aligned}
    \end{equation}
    We may think of $\Ee$ as the only source of randomness in the ``posterior world'' where $\data$ and $\Se$ are fixed, non-random. Theorem \ref{Th: Covariance concentration} applied to $\Ee$ as the difference of the sample covariance of the vectors $\{ Z_j \}_{j=1}^{n_{r,b}}$ and the true covariance $\I_r$ implies that there exist set $\Upsilon$ of the posterior measure 
    \begin{equation}
        \begin{aligned}
            \Pi^{\Wsht}(\Upsilon\,|\,\data) \geq 1 - \frac{1}{n} \text{ for all } \data
            \nonumber
        \end{aligned}
    \end{equation}
    on which
    \begin{equation}
        \begin{aligned}
            \|\Ee\|_{\infty} \lesssim \sqrt{\frac{r + \log(n_{r,b})}{n_{r,b}}} \lesssim \sqrt{\frac{r+\log(n)}{n_{r,b}}}.
            \label{Bound: E}
        \end{aligned}
    \end{equation}
    Without loss of generality we assume that $r/n$ is smaller then some implicit constant (otherwise the claim of the whole theorem is trivial) to guarantee $\| \Ee \|_{\infty} \leq 1/2$.
Now the idea is to get rid of the inversion $(\I_r + \Ee)^{-1}-\I_r$ and to work with $\Ee$ directly. The next lemma allows us to do that. At the same time we get rid of $\Aa$.
\begin{lemma} \label{Lemma: Fourth}
    For the remainder $\mathcal{R}_4$ defined as
    \begin{equation}
        \begin{aligned}
            &\mathcal{R}_4 \eqdef \sqrt{n}\,\Tr\left[ \Aa^{1/2}\D\Aa^{1/2}\cdot \frac{n_{corr}}{n_{r,b}}\left((\I_r + \Ee)^{-1} - \I_r\right)\right] - \sqrt{n_{r,b}}\, \Tr\left[\D\, (-\Ee)\right]
            \nonumber
        \end{aligned}
    \end{equation}
    and
    \begin{equation}
        \begin{aligned}
            \Delta_4 \stackrel{\operatorname{def}}{\asymp}
            \|\D\|_1 \sqrt{r + \log(n)} \left( \wdelta_{n,r} + \sqrt{\frac{r+\log(n)}{n_{r,b}}} + \frac{\|\Gg\|_{\infty}}{n_{r,b}}\right)
            \label{Delta4}
        \end{aligned}
    \end{equation}
    holds
    \begin{equation}
        \begin{aligned}
            \Pi^\Wsht\big(|\mathcal{R}_4| \leq \Delta_4 \;|\;\data \big) \geq 1 - \frac{1}{n}
            \nonumber
        \end{aligned}
    \end{equation}
    with probability $1-1/n$.
\end{lemma}
\begin{proof}
    Firstly, note that $\sqrt{n}\cdot n_{corr}/n_{r,b} = \sqrt{n_{r,b}}$, so the scalers coincide.
    
    Further, since $\| \Ee \|_{\infty} \leq 1/2 < 1$, we can write the Neumann series
    \begin{equation}
        \begin{aligned}
            (\I_r + \Ee)^{-1} = \sum\limits_{k=0}^{\infty} (-1)^k \Ee^k.
            \nonumber
        \end{aligned}
    \end{equation}
    Thus, we have
    \begin{equation}
        \begin{aligned}
            \| (\I_r + \Ee)^{-1} - \I_r - (-\Ee) \|_{\infty} &= 
            \left\| \sum\limits_{k=2}^{\infty} (-1)^k \Ee^k \right\|_{\infty} \leq
            \sum\limits_{k=2}^{\infty} \|\Ee\|_{\infty}^k \\
            &=
            \frac{\|\Ee\|_{\infty}^2}{1-\|\Ee\|_{\infty}} \leq 2\|\Ee\|_{\infty}^2.
            \nonumber
        \end{aligned}
    \end{equation}
    Now we readily bound
    \begin{equation}
        \begin{aligned}
            &\frac{|\mathcal{R}_4|}{\sqrt{n_{r,b}}} =
            \left|\Tr\left[ \Aa^{1/2}\D\Aa^{1/2} \,\left((\I_r + \Ee)^{-1} - \I_r\right)\right] -  \Tr\left[\D\, (-\Ee)\right] \right|\\
            &=
            \Big|\Tr\left[ \Aa^{1/2}\D\Aa^{1/2} \,\left((\I_r + \Ee)^{-1} - \I_r + \Ee\right)\right] + \Tr\left[(\Aa^{1/2}\D\Aa^{1/2} - \D)\, (-\Ee)\right] \Big| \\
            &\leq  \|\Aa^{1/2}\D\Aa^{1/2}\|_1 \| (\I_r + \Ee)^{-1} - \I_r - (-\Ee) \|_{\infty} + \|\Aa^{1/2}\D\Aa^{1/2} - \D\|_1\|\Ee\|_{\infty} \\
            &\lesssim \|\D\|_1 \|\Aa\|_{\infty} \| \Ee \|_{\infty}^2
            + \| \Aa^{1/2}\D \Aa^{1/2} - \Aa^{1/2}\D + \Aa^{1/2} \D - \D \|_1\|\Ee\|_{\infty}\\
            &\leq \|\D\|_1 \|\Aa\|_{\infty} \| \Ee \|_{\infty}^2 + 
            \|\D\|_1 \|\Aa^{1/2}-\I_r\|_{\infty} (\| \Aa\|_{\infty} + 1) \| \Ee \|_{\infty}.
            \nonumber
        \end{aligned}
    \end{equation}
    To finish the proof, we apply Lemma \ref{Lemma: A} and bound (\ref{Bound: E}), and omit higher order terms.
\end{proof}
The final step is Gaussian approximation: we apply the classical Berry-Esseen result to show approximate normality of $\sqrt{n_{r,b}} \Tr\left[\D\, (-\Ee)\right]$ in the next lemma.
\begin{lemma} \label{Lemma: Fifth}
Let $\zeta \sim \mathcal{N}(0, \;2\,\| \D \|_2^2)$. Then
    \begin{equation}
        \begin{aligned}
            &\sup\limits_{x\in\R} 
            \left|
            \Pi^{\Wsht}\left(\sqrt{n_{r,b}}\, \Tr\left[\D\, (-\Ee)\right] \leq x\,\big| \data \right) - \Pro(\zeta \leq x)
            \right| \leq \\
            &\hspace{3cm} \leq \Delta_5 \stackrel{\operatorname{def}}{\asymp} \frac{1}{\sqrt{n}}
            \label{Delta5}
        \end{aligned}
    \end{equation}
with probability $1$.
\end{lemma}
\begin{proof}
    Define the i.i.d. random variables
    \begin{equation}
        \begin{aligned}
            \eta_j \eqdef -Z_j^{\top}\D Z_j = -\sum\limits_{i=1}^r d_i z_{ji}^2, \;\;\;j = 1,\ldots,n_{r,b},
            \nonumber
        \end{aligned}
    \end{equation}
    where $z_{ji}$ is $i$-th component of vector $Z_j$ and $d_i$ is $(i,i)$-th entry of $\D$. Note that since $Z_j$'s are independent standard Gaussian random vectors, all $z_{ji}$'s are independent. 
    Observe that
    \begin{equation}
        \begin{aligned}
            \sqrt{n_{r,b}}\, \Tr\left[\D\, (-\Ee)\right] = \frac{1}{\sqrt{n_{r,b}}} \sum\limits_{j=1}^{n_{r,b}} (\eta_j - \E[\eta_j]).
            \label{Sum}
        \end{aligned}
    \end{equation}
    Let us compute the variance $\sigma^2$ of $\eta_1$:
    \begin{equation}
        \begin{aligned}
            &\sigma^2 \eqdef \Var[\eta_1] = 
            \Var\left[-\sum\limits_{i=1}^r d_i z_{1i}^2\right]
            = \sum\limits_{i=1}^r d_i^2 \, \E\left[(z_{1i}^2 - 1)^2\right] = 2\,\|\D\|_2^2.
        \nonumber
        \end{aligned}
    \end{equation}
    To apply Gaussian approximation to the right-hand side of (\ref{Sum}), we also need to bound the third absolute central moment $\rho$ of $\eta_1$:
    \begin{equation}
        \begin{aligned}
            \rho &= \E\left[ |\eta_1 - \E[\eta_1]|^3\right]
            = \E\left[ \left| \sum\limits_{i=1}^r d_i(z_{1i}^2 - 1)\right|^3\right] \leq \left( \E\left[ \left( \sum\limits_{i=1}^r d_i(z_{1i}^2 - 1) \right)^6\right]\right)^{1/2} \\
            &=\Bigg( 
            \sum\limits_{i=1}^r d_i^6 \, \E\left[(z_{1i}^2-1)^6\right] +
            \sum\limits_{i\neq k} \,d_i^4 \,d_k^2\, \E\left[(z_{1i}^2-1)^4\right]\,\E\left[(z_{1k}^2-1)^2\right] +\\
            &\hspace{1cm}+\sum\limits_{i\neq k} \,d_i^3 \,d_k^3\, \E\left[(z_{1i}^2-1)^3\right]\,\E\left[(z_{1k}^2-1)^3\right] +\\
            &\hspace{1cm}+\sum\limits_{i\neq k, k\neq m, m\neq i} \,d_i^2 \,d_k^2\,d_m^2 \E\left[(z_{1i}^2-1)^2\right]\,\E\left[(z_{1k}^2-1)^2\right]\,\E\left[(z_{1m}^2-1)^2\right]
            \Bigg)^{1/2}.
            \nonumber
        \end{aligned}
    \end{equation}
    When we expanded the sixth power in the display above, we took into account that the other terms disappear due to independence after expectation acts on them. 
    Hence, up to a multiplicative constant we get the upper bound
    \begin{equation}
        \begin{aligned}
            \rho &\lesssim \Bigg( 
            \sum\limits_{i=1}^r d_i^6 +
            \sum\limits_{i\neq k} \,d_i^4 \,d_k^2 +\sum\limits_{i\neq k} \,d_i^3 \,d_k^3+\sum\limits_{i\neq k, k\neq m, m\neq i} \,d_i^2 \,d_k^2\,d_m^2
            \Bigg)^{1/2} \\
            &\leq  \Bigg( 
            \| \D \|_{\infty}^4 \| \D \|_2^2
            + 2\| \D \|_{\infty}^2 \| \D \|_2^4
            + \| \D \|_2^6
            \Bigg)^{1/2} \leq 2\,\| \D\|_2^3.
            \nonumber
        \end{aligned}
    \end{equation}
    The Berry-Esseen Theorem claims that the Kolmogorov distance can be bounded by $\rho/\sigma^3\sqrt{n_{r,b}}$ up to some constant factor. This brings us to the desired statement of the lemma.
\end{proof}

Finally, we put together all the obtained bounds. 
For \( \Delta_1, \Delta_2, \Delta_3 \) from (\ref{Delta1}), (\ref{Delta2}), (\ref{Delta3}), respectively, we write
\begin{equation}
    \begin{aligned}
	&\PpostbigW{ \sqrt{n}\left( \phi(\Su) - \phi(\Se) \right) \leq x } \leq \\
	& \hspace{0.5cm} \leq 
	\PpostbigW{ \sqrt{n} \Tr\left[\Aa^{1/2}\D\Aa^{1/2} \left(\Bb - \frac{n_{corr}}{n_{r,b}}\,\I_r\right)\right] \leq x + \Delta_1 + \Delta_2 + \Delta_3 } \\
	& \hspace{0.5cm} +
	\PpostbigW{ -\sqrt{n}\,\res(\Su,\St) \leq -\Delta_1}
	+
	\PpostbigW{ -\sqrt{n}\,\res(\Se,\St) \leq -\Delta_2} \\
	& \qquad +
	\PpostbigW{  \mathcal{R}_3 \leq -\Delta_3}.
	\nonumber
	\end{aligned}
\end{equation}
The second term in the right-hand side is at most $1/n$ with probability \(1 - 1/n\) due to (\ref{Prob_Delta1}).
The third term is zero with probability $1 - 1/n$ according to (\ref{Prob_Delta2}).
Lemma \ref{Lemma: Third} implies that the fourth term disappears with probability $1 - 1/n$.
Thus, with probability $1 - 3/n$ we have
\begin{equation}
    \begin{aligned}
	&\PpostbigW{ \sqrt{n}\left( \phi(\Su) - \phi(\Se) \right) \leq x } \leq \\
	& \hspace{0.5cm} \leq 
	\PpostbigW{ \sqrt{n} \Tr\left[\Aa^{1/2}\D\Aa^{1/2} \left(\Bb - \frac{n_{corr}}{n_{r,b}}\,\I_r\right)\right] \leq x + \Delta_1 + \Delta_2 + \Delta_3 } \\
	& \hspace{1.5cm} + \frac{1}{n} \\
	& \hspace{0.5cm} = 
	\PpostbigW{ \sqrt{n}\,\Tr\left[ \Aa^{1/2}\D\Aa^{1/2}\cdot \frac{n_{corr}}{n_{r,b}}\left((\I_r + \Ee)^{-1} - \I_r\right)\right] \leq x + \Delta_1 + \Delta_2 + \Delta_3}\\
    &\hspace{1.5cm}+ \frac{1}{n},
	\nonumber
	\end{aligned}
\end{equation}
where the last equality is due to (\ref{Distribution}).
Further, with $\Delta_4$ from (\ref{Delta4}) we have
\begin{equation}
    \begin{aligned}
	&\PpostbigW{ \sqrt{n}\left( \phi(\Su) - \phi(\Se) \right) \leq x } \leq \\
	& \hspace{0.5cm} \leq 
	\PpostbigW{ \sqrt{n}\,\Tr\left[ \Aa^{1/2}\D\Aa^{1/2}\cdot \frac{n_{corr}}{n_{r,b}}\left((\I_r + \Ee)^{-1} - \I_r\right)\right] \leq x + \Delta_1 + \Delta_2 + \Delta_3}\\
	& \hspace{1.5cm} + \frac{1}{n} \\
	& \hspace{0.5cm} \leq
	\PpostbigW{ \sqrt{n}\,\Tr\left[ \D(-\Ee) \right] \leq x + \Delta_1 + \Delta_2 + \Delta_3 + \Delta_4} + \PpostbigW{ \mathcal{R}_4 \leq -\Delta_4 }\\
    &\hspace{1.5cm}+ \frac{1}{n}
	\nonumber
	\end{aligned}
\end{equation}
with probability $1 - 3/n$. The second term in the right-hand side is at most $1/n$ with probability $1-1/n$ by Lemma \ref{Lemma: Fourth}, thus
\begin{equation}
    \begin{aligned}
	&\PpostbigW{ \sqrt{n}\left( \phi(\Su) - \phi(\Se) \right) \leq x } \leq \\
	& \hspace{0.5cm} \leq
	\PpostbigW{ \sqrt{n}\,\Tr\left[ \D\,(-\Ee) \right] \leq x + \Delta_1 + \Delta_2 + \Delta_3 + \Delta_4} + \frac{2}{n}
	\nonumber
	\end{aligned}
\end{equation}
with probability $1 - 4/n$.

Subtracting $\Pro(\zeta \leq x)$ with $\zeta \sim \ND(0, 2\| \D \|_2^2)$ from the both sides, inserting additional $\pm\Pro(\zeta \leq x + \Delta_1 + \Delta_2 + \Delta_3 + \Delta_4)$ and taking supremum over $x\in \R$, we obtain
\begin{equation}
    \begin{aligned}
	&\sup\limits_{x\in\R} \left[ \PpostbigW{ \sqrt{n}\left( \phi(\Su) - \phi(\Se) \right) \leq x } - \Pro(\zeta \leq x) \right] \leq \\
	& \hspace{1cm} \leq 
	\sup\limits_{x\in\R} 
	\Bigg[ \PpostbigW{ -\sqrt{n}\, \Tr\left[ \D\,(-\Ee) \right] \leq x + \Delta_1 + \Delta_2 + \Delta_3 + \Delta_4} - \\
	& \hspace{6.3cm}
	- \Pro(\zeta \leq x+\Delta_1 + \Delta_2 + \Delta_3 + \Delta_4)\Bigg] + \\
	&\hspace{2cm} + 
	\sup\limits_{x\in\R}
	\left[ \Pro(\zeta \leq x+\Delta_1 + \Delta_2 + \Delta_3 + \Delta_4) - \Pro(\zeta \leq x)\right]
	+ \frac{2}{n}
    \nonumber
	\end{aligned}
\end{equation}
with probability $1 - 4/n$.
The first term in the right-hand side is bounded by $\Delta_5$ from (\ref{Delta5}) with probability $1$ due to Lemma~\ref{Lemma: Fifth}.
The second term is at most \\
$(\Delta_1 + \Delta_2 + \Delta_3 + \Delta_4)/(\sqrt{2\pi}\,\sqrt{2} \| \D \|_2)$.

Hence,
\begin{equation}
    \begin{aligned}
	&\sup\limits_{x\in\R} \left[ \PpostbigW{ \sqrt{n}\left( \phi(\Su) - \phi(\Se) \right) \leq x } - \Pro(\zeta \leq x) \right] \leq \\
	& \hspace{1cm} \leq 
	\frac{\Delta_1 + \Delta_2 + \Delta_3 + \Delta_4}{\| \D \|_2}+ \Delta_5 + \frac{2}{n}
    \nonumber
	\end{aligned}
\end{equation}
with probability $1 - 4/n$. Similarly, one can show
\begin{equation}
    \begin{aligned}
	&\sup\limits_{x\in\R} \left[ \Pro(\zeta \leq x) - \PpostbigW{ \sqrt{n}\left( \phi(\Su) - \phi(\Se) \right) \leq x } \right]\\
	& \hspace{1cm} \leq 
	\frac{\Delta_1 + \Delta_2 + \Delta_3 + \Delta_4}{\| \D \|_2}+ \Delta_5 + \frac{2}{n}
    \nonumber
	\end{aligned}
\end{equation}
with probability $1 - 4/n$.
Observe that 
\begin{equation}
    \begin{aligned}
    \frac{\|\D\|_1}{\| \D\|_2} \leq \sqrt{r}  \;\text{  and  }\; \| \D\|_2 = \left\| \St^{1/2} \Phib \St^{1/2} \right\|_2.
    \nonumber
	\end{aligned}
\end{equation}
By adjusting the constants, we get the desired statement with probability $1~-~1/n$.

	\appendix
	\section{Auxiliary results} \label{Appendix: A}
		Here we collect some auxiliary statements from the literature.

The next Theorem gathers several results on concentration of sample covariance, which allows us to work in more general framework than just Gaussian distribution of the observations. The same collection of results presented in \cite{Silin}.
\begin{theorem} \label{Th: Covariance concentration}
	Let \( X_1, \ldots, X_n \) be i.i.d. zero-mean random vectors in \( \R^p \). 
	Denote the true covariance matrix as \( \St \eqdef \E \left(X_j X_j^{\T}\right) \) and the sample covariance as \\ \( \Se \eqdef \frac{1}{n}\sum_{j=1}^n X_j X_j^{\T} \).
	Suppose the data and are obtained from:\\
	(i) Gaussian distribution \( \ND(0, \St) \). In this case, define \( \hdelta_{n,p} \) as
	\begin{equation}
		\begin{aligned}
			\hdelta_{n,p} \asymp  \sqrt{\frac{\reff(\St) + \log(n)}{n}} ;
		\nonumber
		\end{aligned}
	\end{equation}
	(ii) sub-Gaussian distribution. In this case, define \( \hdelta_{n,p} \) as
	\begin{equation}
		\begin{aligned}
			\hdelta_{n,p} \asymp \sqrt{\frac{p + \log(n)}{n}};
		\nonumber
		\end{aligned}
	\end{equation}
	(iii) a distribution supported in some centered Euclidean ball of radius \( R \). In this case, define \( \hdelta_{n,p} \) as
	\begin{equation}
		\begin{aligned}
			\hdelta_{n,p} \asymp  \frac{R}{\sqrt{\| \St \|}} \sqrt{\frac{\log(n)}{n}};
		\nonumber
		\end{aligned}
	\end{equation}
	(iv) log-concave probability measure. In this case, define \( \hdelta_{n,p} \) as
	\begin{equation}
		\begin{aligned}
			\hdelta_{n,p} \asymp  \sqrt{\frac{\log^6(n)}{np}} .
		\nonumber
		\end{aligned}
	\end{equation}
	
	Then in all the cases above the following concentration result for \( \Se \) holds with the corresponding \( \hdelta_{n,p} \):
	\begin{equation}
		\begin{aligned}
			\|\Se - \St\|_{\infty} \leq \hdelta_{n,p} \| \St\|_{\infty}
		\nonumber
		\end{aligned}
	\end{equation}
	with probability at least \( 1-\frac{1}{n} \).
\end{theorem}
\begin{proof}
	(i) See \cite{Koltchinskii_CIAMBFSCO}, Corollary 2.
	(ii) This is a well-known simple result presented in a range of papers and lecture notes. See, e.g. \cite{Rigollet}, Theorem 4.6.
	(iii) See \cite{Vershynin_ITTNAAORM}, Corollary 5.52. Usually the radius \( R \) is taken such that \( \frac{R}{\sqrt{\| \St \|}} \asymp \frac{\sqrt{\Tr(\St)}}{\sqrt{\| \St \|}} = \sqrt{\widetilde{r}(\St)} \).
	(iv) See \cite{Adamczak}, Theorem 4.1.
\end{proof}

Throughout the paper we use the properties of the Inverse Wishart distribution presented in the following lemma:
\begin{lemma} \label{Lemma: Properties}
    Let $\W \sim \Wsht_k^{-1}(\Hh,\, s)$. Then the following holds:\\
    (i) \;\;\;\;$\E[\W] = \frac{\Hh}{s-k-1}$;\\
    (ii) \;\;\;$\Q^{\top} \W \Q \sim \Wsht_q^{-1}(\Q^{\top} \Hh \Q,\, s-k+q)$ for any $\Q\in\R^{k\times q}$;\\
    (iii) \;\;$\W^{-1} \stackrel{d}{=} \sum\limits_{j=1}^s Z_j Z_j^{\top}$ for $\{ Z_j \}_{j=1}^s \stackrel{i.i.d.}{\sim} \mathcal{N}_k(0, \Hh)$.
\end{lemma}
\begin{proof}
    We refer to Chapter 3 of \cite{Muirhead} which is a great reference for properties of Wishart and Inverse Wishart distributions.
\end{proof}

	\section{Auxiliary proofs} \label{Appendix: B}
	 	Here we present proofs of Proposition \ref{Contraction} and Proposition \ref{Flatness}. 
\begin{proof}[Proof of Proposition \ref{Contraction}]
We have $\Su \sim \Wsht_p^{-1}(\G,\, p + b -1)$, and the formulae for the prior density $\Pi^{\Wsht}$ and the Gaussian log-likelihood $l_n$ imply the posterior density which corresponds to
\begin{equation}
    \begin{aligned}
        \Su \,|\, \data \sim \Wsht_p^{-1}(n\Se + \G,\, n+p+b-1).
        \nonumber
    \end{aligned}
\end{equation}
Introduce for shortness
\begin{equation}
    \begin{aligned}
        n_{p,b} &\eqdef n+p+b-1,\\
        \Su_{\circ} &\eqdef \frac{n\Se+\G}{n_{p,b}}.
        \nonumber
    \end{aligned}
\end{equation}
Decompose
\begin{equation}
    \begin{aligned}
        \Su - \St &= (\Su - \Su_{\circ}) + (\Su_{\circ} - \St) \\
        &=\Su_{\circ}^{1/2}\, (\Su_{\circ}^{-1/2}\Su\Su_{\circ}^{-1/2}-\I_p) \,\Su_{\circ}^{1/2} + (\Su_{\circ} - \St).
        \nonumber
    \end{aligned}
\end{equation}

Let us bound the first term. Note that by Lemma \ref{Lemma: Properties} (ii)
\begin{equation}
    \begin{aligned}
        \frac{1}{n_{p,b}}\,\Su_{\circ}^{-1/2}\Su\Su_{\circ}^{-1/2} \,|\, \data \sim \Wsht_p^{-1}(\I_p,\, n_{p,b}).
        \nonumber
    \end{aligned}
\end{equation}
So, with 
    \begin{equation}
        \begin{aligned}
            \Ee = \frac{1}{n_{p,b}} \sum\limits_{j=1}^{n_{p,b}} Z_j Z_j^{\top} - \I_p,\;\;\{ Z_j \}_{j=1}^{n_{p,b}}\,|\,\data \stackrel{i.i.d.}{\sim} \mathcal{N}_p(0, \I_p),
    \nonumber
        \end{aligned}
    \end{equation}
    Lemma \ref{Lemma: Properties} (iii) yields
    \begin{equation}
        \begin{aligned}
            \Su_{\circ}^{-1/2}\Su\Su_{\circ}^{-1/2} - \I_p\;\Big|\;\data \eqdist (\I_p + \Ee)^{-1} - \I_p.
            \nonumber
        \end{aligned}
    \end{equation}
    By Theorem \ref{Th: Covariance concentration}, there exist a set $\Upsilon$ of posterior measure
    \begin{equation}
        \begin{aligned}
            \Pi^{\Wsht} (\Upsilon \;|\;\data) \geq 1-\frac{1}{n}\;\;\text{for all } \data
            \nonumber
        \end{aligned}
    \end{equation}
    on which
    \begin{equation}
        \begin{aligned}
            \| \Ee \|_{\infty} \lesssim \sqrt{\frac{p+\log(n_{p,b})}{n_{p,b}}} \leq \sqrt{\frac{p+\log(n)}{n}}.
            \nonumber
        \end{aligned}
    \end{equation}
    Now we use the assumption that $p/n$ is small enough to provide $\| \Ee \|_{\infty} \leq 1/2$. Then, using Neumann series, on $\Upsilon$ holds
    \begin{equation}
        \begin{aligned}
            \| (\I_p + \Ee)^{-1} - \I_p \|_{\infty} &= 
            \left\| \sum\limits_{k=1}^{\infty} (-1)^k \Ee^k \right\|_{\infty} \leq
            \sum\limits_{k=1}^{\infty} \|\Ee\|_{\infty}^k \\
            &=
            \frac{\|\Ee\|_{\infty}}{1-\|\Ee\|_{\infty}} \leq 2\,\|\Ee\|_{\infty} \lesssim \sqrt{\frac{p+\log(n)}{n}}.
            \nonumber
        \end{aligned}
    \end{equation}
    Hence,
    \begin{equation}
        \begin{aligned}
            \| \Su - \Su_{\circ}\|_{\infty} &=
            \| \Su_{\circ}^{1/2}\, (\Su_{\circ}^{-1/2}\Su\Su_{\circ}^{-1/2}-\I_p) \,\Su_{\circ}^{1/2}\|_{\infty} \\
            &\leq
            \| \Su_{\circ} \|_{\infty} \|\Su_{\circ}^{-1/2}\Su\Su_{\circ}^{-1/2}-\I_p\|_{\infty}\\
             &\lesssim (\| \St \|_{\infty} + \| \Su_{\circ} - \St \|_{\infty})\, \sqrt{\frac{p+\log(n)}{n}}
            \nonumber
        \end{aligned}
    \end{equation}
    on a set of posterior measure $1-1/n$.

Now we bound the second term:
    \begin{equation}
        \begin{aligned}
            \| \Su_{\circ} - \St \|_{\infty} &=
            \left\| \frac{n\Se+\G}{n_{p,b}} - \St \right\|_{\infty} \\
            &\leq \| \Se - \St\|_{\infty} + \frac{p+b}{n_{p,b}}\,\|\St\|_{\infty} + \frac{\|\G\|_{\infty}}{n_{p,b}}\\
            &\leq \hdelta_{n,p} \,\| \St\|_{\infty} + \frac{p+b}{n}\,\|\St\|_{\infty} + \frac{\|\G\|_{\infty}}{n},
            \nonumber
        \end{aligned}
    \end{equation}
where the last inequality holds with probability $1-1/n$ due to the covariance concentration condition (\ref{Eq: hat_delta_n}). 
Putting these bounds together and omitting higher-order terms, we deduce that the posterior measure of the event
\begin{EQA}
	\Pi^\Wsht\left( \| \Su - \St \|_{\infty} 
	\lesssim 
	\| \St \|_{\infty} \left( \sqrt{\frac{\log(n)+p}{n}} + \hdelta_n + \frac{\|\G\|_{\infty}}{n\,\|\St\|_{\infty}} \right) \,\Bigg|\, \data\right) \geq 1-\frac{1}{n}
\end{EQA}
with probability $1 -  1/n$.
\end{proof}

\begin{proof}[Proof of Proposition \ref{Flatness}]
	The log-density of the Inverse Wishart prior is
	\begin{EQA}
	    	\log{\Pi^{\Wsht}(\Su)} 
	        & = & 
	        -\frac{2p+b}{2} \log{\det(\Su)} - \frac{1}{2} \Tr(\G\Su^{-1}) + C,
	\end{EQA}
	where $C$ is some constant related to normalization. 
	For an arbitrary $\Su \in \vic(\delta)$ we have
	\begin{equation}
	    \begin{aligned}
	    	&\left| \frac{\Pi^{\Wsht}(\Su)}{\Pi^{\Wsht}(\St)} - 1\right| 
	    	= 
	    	\left| \exp{ \left\{ \log{\Pi^{\Wsht}(\Su)} - \log{\Pi^{\Wsht}(\St)} \right\}} - 1\right| \\
	    	& = 
	    	\left| \exp{ \left\{ -\frac{2p+b}{2} \left(\log{\det(\Su)} - \log{\det(\St)} \right) - \frac{1}{2} \Tr\left[\G\left(\Su^{-1}-\St^{-1}\right)\right] \right\}} - 1 \right| \\
	    	& \lesssim \
	    	\exp{ \left\{ (p+b/2)\,\left| \log{\det(\Su)} - \log{\det(\St)} \right| + \frac{1}{2}\left| \Tr\left[\G\left(\Su^{-1}-\St^{-1}\right)\right] \right| \right\}} - 1.
	    \nonumber
	    \end{aligned}
	\end{equation}
	Notice that
	\begin{EQA} 
	    	\left| \log{\det(\Su)} - \log{\det(\St)} \right| 
	        & = & 
	        \log{\det({\St}^{-1/2} \;\Su\; {\St}^{-1/2})} 
	        = \log{\det(\I_p + \boldsymbol{\Delta})},
	\end{EQA}
	where $ \boldsymbol{\Delta} \eqdef \St^{-1/2} (\Su-\St) \St^{-1/2}$ and $\| \boldsymbol{\Delta} \|_{\infty} \leq  \| {\St^{-1}} \|_{\infty} \,\| \St \|_{\infty}\, \delta \leq 1/2$ by our assumption. Then
	\begin{EQA} 
	    	\log{\det(\I_p + \boldsymbol{\Delta})} 
	    	= 
	    	\sum\limits_{j=1}^p \log{(1+\lambda_j(\boldsymbol{\Delta}))} 
	        & = & 
	    	-\sum\limits_{j=1}^p \sum\limits_{s=1}^{\infty} \frac{\lambda_j^s(\boldsymbol{\Delta})}{s} 
	    	= 
	    	-\sum\limits_{s=1}^{\infty} \frac{\Tr(\boldsymbol{\Delta}^s)}{s}.
	\end{EQA}
	Hence, taking into account $ p\,\delta\,\| {\St^{-1}} \|_{\infty}\,\| \St \|_{\infty} \leq 1/2$, we get the bound
	\begin{equation}
	    \begin{aligned}
	    	\left| \log{ \det(\I_p + \boldsymbol{\Delta})} \right| 
	        & \leq 
	    	\sum\limits_{s=1}^{\infty} \| {\St^{-1}} \|_1^s \, \| \Su-\St\|_{\infty}^s 
	    	= 
	    	\frac{\| {\St^{-1}} \|_1\,\| \St \|_{\infty} \, \delta}{1 - \| {\St^{-1}} \|_1\,\| \St \|_{\infty} \, \delta} \\
	    	& \leq 2\,p\,\delta\,\| {\St^{-1}} \|_{\infty}\,\| \St \|_{\infty}.
	    	\nonumber
	    \end{aligned}
	\end{equation}
	Besides,  for $ \Su \in \vic(\delta) $ we have
	\begin{EQA}
		\| \Su^{-1} - {\St}^{-1} \|_{\infty}
		& = &
		\| {\St}^{-1} (\Su - \St) \Su^{-1} \|_{\infty} 
		\leq
		\delta \,\| \St \|_{\infty}\, \| {\St}^{-1} \|_{\infty} \,\| \Su^{-1} \|_{\infty}\\
		& \leq &
		\delta \,\| \St \|_{\infty}\, \| {\St}^{-1} \|_{\infty} \,(\| {\St}^{-1} \|_{\infty} +\| \Su^{-1} - {\St}^{-1} \|_{\infty}).
	\end{EQA}
	Therefore,
	\begin{EQA}
		\| \Su^{-1} - {\St}^{-1} \|_{\infty}
		& \leq &
		\frac{\delta \,\| \St \|_{\infty}\, \| {\St^{-1}} \|_{\infty}^2}{1 - \delta\,\| \St \|_{\infty}\, \| {\St^{-1}} \|_{\infty}}
		\leq
		2\,\delta\,\| \St \|_{\infty} \, \| {\St^{-1}} \|_{\infty}^2.
	\end{EQA}
	Hence,
	\begin{EQA}
		\left| \Tr\left[\G \left( \Su^{-1} - {\St}^{-1} \right)\right] \right| 
		& \leq & 
		2\,\delta \, \| \G \|_1 \| \St \|_{\infty} \| {\St^{-1}} \|_{\infty}^2.
	\end{EQA}
	Finally,
	\begin{EQA} 
	    	\fl^{\Wsht}(\delta) = \sup\limits_{\Su \in \vic{\delta}} \left| \frac{\Pi^{\Wsht}(\Su)}{\Pi^{\Wsht}(\St)} - 1\right|
	        & \lesssim &
	        \delta \,\kappa(\St) \left\{ p^2 + \| \G \|_1 \| {\St^{-1}} \|_{\infty} \right\}.
	\end{EQA}
\end{proof}

	\section*{Acknowledgements}
	The author is grateful to Vladimir Spokoiny for valuable discussions and to Chao Gao for helpful comments.

	\bibliographystyle{apalike}
	
\end{document}